\documentclass[11pt]{article} 
\usepackage{amsmath, amscd}
\usepackage{amssymb, amsthm, latexsym} 

\usepackage[all]{xy}
\usepackage{epic,eepic}

\CompileMatrices

\newcommand{\PP}{{\mathbb P}}
\newcommand{\QQ}{{\mathbb Q}}
\newcommand{\ZZ}{\mathbb Z}
\newcommand{\Gr}{{\mathbb G}}

\newcommand{\Mg}{\overline{{\cal M}}_g}
\newcommand{\Mgi}{\overline{{\cal M}}_{g,1}}
\newcommand{\Mgn}{\overline{{\cal M}}_{g,n}}

\newcommand{\Moi}{\overline{{\cal M}}_{0,1}}
\newcommand{\Grk}{{\Gr(E,k)}}
\newcommand{\Grpk}{{\Gr(E_p,k)}}

\newcommand{\Tvert}{{T_{\mbox{\scriptsize\sc V}}}}
\newcommand{\Hom}{{\mbox{Hom}}}
\newcommand{\Ext}{{\mbox{Ext}}}

\newcommand{\Quot}{\mbox{Quot}_{k,d}(E)}
\newcommand{\Cfam}{{\mathcal{C}}}
\newcommand{\Qfam}{{\mathcal{Q}}}
\newcommand{\Vbndle}{{\mathcal{V}}}
\newcommand{\Sbndle}{{\mathcal{S}}}
\newcommand{\Fbndle}{{\mathcal{F}}}
\newcommand{\Osh}{{\cal{O}}}
\newcommand{\Id}{{\cal{I}}}
\newcommand{\pr}{\mbox{pr}}
\newcommand{\doneproof}{{\rule{2mm}{2mm}}}

\newlength{\sectiontitlewidth}
\newcommand{\secheading}[1]{{\parbox[t]{\sectiontitlewidth}{\Large #1}}}
\newcommand{\mysection}[1]{{\Large\bf{\refstepcounter{section}}\arabic{section}. 
\secheading{#1}}\addcontentsline{toc}{subsection}{\protect{\arabic{section}. #1}}}
\newcounter{example}[section]
\newcommand{\example}{\refstepcounter{example}{\bf Example \arabic{section}.\arabic{example}: }}

\newtheorem{theorem}{Theorem}[section]
\newtheorem{lemma}{Lemma}[section]

\newtheorem{corollary}[theorem]{Corollary}
\newtheorem{proposition}[theorem]{Proposition}

\begin{document}

\title{\bf Stable maps and Quot schemes}
\author{Mihnea Popa {\rm and} Mike Roth}
\date{}

\maketitle 

\tableofcontents

\setlength{\parskip}{.1 in} 

\bigskip
\noindent
\mysection{Introduction}

\medskip
\noindent
Let $C$ be a smooth projective curve of genus $g$ over an algebraically closed
field of arbitrary characteristic, and $E$ a vector bundle on $C$
of rank $r$. The main theme of this paper is the relationship 
between two different 
compactifications of the space of vector bundle quotients of $E$ of fixed rank and degree. 
The first is 
Grothendieck's Quot scheme ${\rm Quot}_{k,d}(E)$ of coherent quotients  
of rank $k$ and degree $d$, while the second is a compactification $\Mg(\Grk, \beta_d)$
via stable maps into the relative 
Grassmannian $\Grk$ of $k$-dimensional quotient spaces of $E$.
Quot schemes have played a central role in the 
``classical'' study of vector bundles on curves, in particular
for constructing moduli spaces (see e.g. \cite{Seshadri}) and for understanding their 
geometry (see \cite{Drezet}). In some general sense,
any construction or computation in this area depends on geometric or 
dimensional properties of (subsets of) Quot schemes. The main product of the study 
of the relationship with the corresponding space of stable maps will be a 
surprisingly good understanding of these properties, as we will explain below.

\medskip
\noindent
Before the availability of stable maps, 
Quot schemes have also been used in Gromov-Witten theory as a way to 
compactify the
space of maps into Grassmannians and thereby compute intersection numbers 
(most notably in \cite{Bertram1}, \cite{Bertram2}, \cite{BDW}
and subsequent extensions \cite{Chen}, \cite{CF1}, \cite{CF2}, \cite{Kim}).
The introduction by Kontsevich of the moduli space of stable maps 
now gives an alternate setting for these calculations; however
in this paper we follow in some respect an orthogonal direction by showing 
how an understanding of the stable map space leads to good results about 
Quot schemes. For instance, we establish an essentially optimal upper bound on the dimension 
of $\Mg(\Grk, \beta_d)$, which gives as an 
immediate consequence an identical upper bound on the dimension of ${\rm Quot}_{k,d}(E)$.
Based on this, we prove what should be seen as the central result of the paper 
-- apparently not even conjectured before -- namely 
that for every vector bundle $E$, $\Quot$ is irreducible and generically smooth
for $d$ sufficiently large, and we precisely describe all the vector bundles for which 
the same thing holds for $\Mg(\Grk, \beta_d)$. (We also make these results effective in
the case of stable bundles.)
We show that, despite their similar behaviour, there are in general no natural 
morphisms between these two compactifications. 
Finally, as a second application of the dimension bound, we obtain new cases 
of a conjecture on the effective base point freeness of linear 
series on moduli spaces of vector bundles on $C$.  

\bigskip
\noindent
We now turn to a more detailed discussion of the results of the paper.

\bigskip
\noindent
In recent literature there has been some interest in
finding practical dimension bounds for families of sheaves or vector bundles on
curves, 
most notably for Quot schemes (see e.g. \cite{Brambila}, \cite{Holla}, \cite{Popa} 
or \cite{Teixidor}). 
Our first goal in this paper is to give a good answer to this question for an 
arbitrary vector bundle $E$.
We do this by considering the second compactification of the space of 
quotient 
bundles of $E$ of rank $k$ and degree $d$, namely the Kontsevich moduli space 
$\Mg(\Grk, \beta_d)$ of stable maps of arithmetic genus $g$ to the relative 
Grassmann bundle $\Grk$ over 
$C$. Here $\beta_d \in H_{2}(\Grk, \ZZ)$ denotes the 
homology class of a section 
corresponding to such a quotient of $E$. Set
$$d_k = d_k (E):=\underset{{\rm rk}\,F=k}{{\rm min}}\left\{{\rm deg}\,F~|\,\,\,\mbox{$F$ 
a quotient of $E$}\right\},$$  
i.e. the minimal degree of a quotient bundle of rank $k$. Our first main result is:

\medskip
\noindent
{\bf Theorem.}
$$ \dim\Mg(\Grk,\beta_d) \le k(r-k) + (d-d_k)r, ~for~ all~ d\ge d_k.$$  

\medskip
\noindent
The advantage of using the stable map space comes from the possibility of 
applying methods characteristic to the study of families of curves.
The idea is that once the space of maps is large enough, the Bend and Break Lemma of Mori 
ensures the existence of degenerate curves related to quotient bundles
of smaller degree, and this gives an inductive bound.
We expect that this approach might be useful in other situations involving 
fibrations over curves.

\medskip
\noindent
Going back to Quot schemes, the theorem above implies almost immediately the 
analogous statement:

\medskip
\noindent
{\bf Theorem.}
$$\dim \Quot \leq k(r-k)+ (d-d_{k})r, ~for~ all~ d\ge d_k.$$

\medskip
\noindent
The case $d=d_k$ of this result was first proved by Mukai-Sakai \cite{Mukai}
and was generalized to the setting of $G$-bundles 
in \cite{Holla}. The general statement improves a similar result 
found in \cite{Popa} by different methods. It is worth emphasizing that 
while this is the statement which will be used in subsequent results, 
it is also here that it is particularly important to first go through the study 
of the space of stable maps: we don't know how to prove the estimate by 
looking only at the Quot scheme.

\medskip
\noindent
We next look at some examples showing certain aspects of the behaviour 
of $\Mg (\Grk, \beta_d)$ and $\Quot$. 
These show that the bounds above are essentially 
optimal if we are not willing to impose extra conditions on the bundle $E$
(and close to optimal even if we are). 
An illustration of this is the fact that the difference between 
our upper bound and the usual deformation theoretic 
lower bound $rd -ke - k(r-k)(g-1)$ gives:
$$rd_k - ke \leq k(r-k)g,$$
which is precisely a theorem of Lange and Mukai-Sakai generalizing a well known result of
Segre and Nagata on minimal sections of ruled surfaces (see \S5). 

\bigskip
\noindent
The dimension estimates also allow us to understand the 
component structure of 
$\Mg (\Grk, \beta_d)$ and $\Quot$. For an arbitrary vector bundle, the picture turns out to 
be as nice as one could possibly hope for.

\medskip
\noindent
{\bf Theorem.} 
\emph{For all large $d$ there is a unique component of  
$\Mg(\Grk,\beta_d)$ whose generic point corresponds to a smooth section 
(i.e. a vector bundle quotient);
this component is of dimension $rd-ke-k(r-k)(g-1)$}.

\medskip
\noindent
Our main result gives an even stronger statement in the case of Quot schemes. Note that 
it is easy to (and we will) construct examples where the Quot scheme is reducible 
for low degree.

\medskip
\noindent
{\bf Theorem.}
\emph{For any vector bundle $E$ on $C$, there is an integer $d_Q = d_Q(E,k)$ such that
for all $d\geq d_Q$, ${\rm Quot}_{k,d}(E)$ is irreducible. For any such $d$,
${\rm Quot}_{k,d}(E)$ is generically smooth, has dimension $rd-ke-k(r-k)(g-1)$, 
and its generic point corresponds to a vector bundle quotient}.

\medskip
\noindent
In the case when $E$ is the trivial bundle (and over the complex numbers),
this result was proved by Bertram-Daskalopoulos-Wentworth (in \cite{BDW} Theorem 4.28), 
and it was used to 
produce an effective way of computing Gromov-Witten invariants on the space of 
holomorphic maps into Grassmannians.

\medskip
\noindent
The proofs of the two theorems rely on the estimates 
above and on some vanishing results of independent interest, involving stable quotients
of arbitrary vector bundles. These are quite elementary, but seem to have been 
overlooked in the literature. Note also that an immediate consequence of the theorem 
is that every vector bundle can be written as an extension of generic stable 
bundles (see Corollary 6.3).

\medskip
\noindent
The next natural question to ask is if it is possible to give 
a precise description of all the bundles for which the spaces  
$\Mg(\Grk,\beta_d)$ are irreducible for all $k$ and all large $d$.
It turns out that, in contrast with the Quot scheme situation,
they are all stable bundles of a very special kind. 
It is shown in 
\cite{Teixidor} that there exists a nonempty open subset of the moduli space $U_C^s(r,e)$
of stable bundles of rank $r$ and degree $e$ whose elements $E$ satisfy the 
properties that $d_k(E)$ is the smallest value making the expression $d_k r -ke - k(r-k)(g-1)$ 
nonnegative, and that ${\rm dim~}{\rm Quot}_{k,d}(E) = rd - ke - k(r-k)(g-1)$ (this 
also applies to $\Mg(\Grk,\beta_d)$). This subset is contained in the better known 
open subset of Lange generic stable bundles (for precise details see Example 5.\ref{generic_stable_bundle}).
Our result is: 

\medskip
\noindent
{\bf Theorem.}
\emph{For any vector bundle $E$ the spaces $\Mg(\Grk,\beta_d)$ are 
connected for all large $d$. 
They are irreducible for all $k$ and
all large $d$ if and only if $E$ is stable bundle, generic in the above sense.} 

\medskip
\noindent
For stable bundles the invariants involved in the proof can be controlled, and in fact a $d$ 
which guarantees irreducibility can be expressed effectively. This is the content of Theorem 
\ref{effective_bound}. A number of other effective statements are given at the end of \S6.

\bigskip
\noindent
Our methods show that in most situations it is more convenient to work with stable 
maps with reducible domain rather than quotients with torsion. In the hope for a 
best possible picture, we can ask whether there is always a morphism $\Mg(\Grk,\beta_d)\rightarrow
{\rm Quot}_{k,d}(E)$ extending the identification on the locus corresponding to 
vector bundle quotients (or, equivalently, smooth sections). We prove that in general this 
question has a negative answer, and so a relationship between the two compactifications, if any,
has to be more indirect. More precisely:

\medskip
\noindent
{\bf Theorem.}\smallskip\newline
(a) \emph{If $k=(r-1)$ then there is a surjective morphism from $\Mg(\Grk,\beta_d)$ to
${\rm Quot}_{k,d}(E)$ which extends the map on the locus where the domain curve is smooth.
Such a morphism also exists for any $k$ if $d = d_k$ or $d = d_k +1$.}
\smallskip\newline
(b) \emph{If $k\ne (r-1)$ then in general there is no morphism from 
$\Mg(\Grk,\beta_d)$
to ${\rm Quot}_{k,d}(E)$ extending the map on smooth curves.}

\bigskip
\noindent
Our final application of the dimension bound discussed above 
is to derive bounds on 
base point freeness for linear series on the moduli spaces of vector bundles
on curves. The fact that uniform 
(i.e., independent of special conditions on $E$) estimates
of this type produce base point freeness bounds is explained in \cite{Popa} \S4.
If $SU_C(r)$ is the moduli space of semistable bundles of rank $r$ with trivial
determinant, and $\mathcal{L}$ is the determinant line bundle, then (cf. \S8.):

\medskip
\noindent
{\bf Theorem. } 
\emph{The linear series $|\mathcal{L}^{p}|$ is base point free 
on $SU_{C}(r)$ for $p\geq [\frac{r^{2}}{4}]$}.

\noindent
This can then be used to show a few
new cases of the Conjecture 5.5 in \cite{Popa} on optimal effective base 
point freeness on the moduli space $U_{C}(r,0)$ of semistable bundles of 
rank $r$ and degree $0$, generalizing the classical case of Jacobians. 
Namely, if $\Theta_{N}$ is the generalized theta divisor on $U_{C}(r,0)$ 
associated to a line bundle $N\in {\rm Pic}^{g-1}(X)$, then we have:

\medskip
\noindent
{\bf Corollary.}
\emph{The linear series 
$|k\Theta_{N}|$ is base point free on $U_{C}(r,0)$ for $k\geq r+1$ if $r\leq 5$.}

\bigskip
\noindent
The paper is organized as follows. In \S2 we describe the basic setup for moduli spaces of stable 
maps in the case of relative Grassmannians. The dimension estimate for $\Mg(\Grk, \beta_d)$ 
is proved in \S3. 
In \S4 we show the analogous bound for Quot schemes. 
In \S5 we give a series of examples leading to a discussion of the optimality 
of the bound. The study of the irreducibility of $\Mg(\Grk, \beta_d)$ and $\Quot$ appears 
in \S6, and in \S7 we show that there are in general no natural morphisms between the 
two compactifications, except in the case $k = (r-1)$. We conclude in \S8 by using the 
dimension estimate to prove the result on base point freeness on moduli spaces. 

\medskip
\noindent
{\bf Acknowledgments.} We would especially like to thank W. Fulton and
R. Lazarsfeld for valuable suggestions and encouragement; W. Fulton  
pointed out a number of references which allowed us to improve the 
results in \S6 and to avoid an error.
We also thank M.S. Narasimhan and M. Teixidor i Bigas for some very useful conversations.  
The second author would also like to express his gratitude to
the Max-Plank institute for a wonderful working environment.
Finally, we are grateful to the referee for some very useful suggestions 
which led to an improvement of the exposition.

\bigskip
\noindent
\mysection{Basic Setup}

\medskip
\noindent
Let $C$ be a smooth projective curve of genus $g$ .
Suppose that $E$ is a vector bundle of rank $r$ and degree $e$ on $C$. For any $k$, $1\le k < r$,  let 
$\Grk$ be the relative Grassmannian of $k$-dimensional quotient spaces of $E$; the fibre of $\Grk$ over a point $p$ of $C$
is the Grassmannian $\Grpk$ of $k$-dimensional quotient spaces of $E_p$.  Let $\pi$ be the projection
$$\pi:\Grk \longrightarrow C.$$

\medskip
\noindent
The quotient bundles $F$ of $E$ of rank $k$ are in one to one correspondence with sections $\sigma_{F}$ of $\pi$.  We want to
compactify the space of quotient bundles by compactifying the space of sections.  One thing we will do to begin with is 
see how the data of the section $\sigma_{F}$ determines the only other piece of numerical data about $F$, its {degree}.

\medskip
\noindent
Let $\Tvert$ be the vertical tangent bundle relative to the map $\pi$.  Also, 
for any quotient bundle $F$, let $S_F$ be the induced subbundle of $E$:
$$0\longrightarrow S_{F}\longrightarrow E\longrightarrow F\longrightarrow 0.$$ 
The pullback $\sigma_{F}^{*}\,\Tvert$ is the bundle $\mathcal{H}\emph{om}(S_F,F)$ on $C$. 
If $F$ has degree $d$, then $\sigma_{F}^{*}\,T_{V}\cong S_{F}^{*}\otimes F$ has degree $(rd - ke)$.
This shows that:
$$ d = \mbox{deg}(F) = \frac{\mbox{deg}(\sigma_{F}^{*}\,\Tvert)+ke}{r}.$$

\medskip
\noindent
Alternatively, if $\beta_F$ is the
homology class $\sigma_{F*}[C]\in H_2(\Grk,{\mathbb Z})$ 
of the section, then we can also write this as:
$$ d = \mbox{deg}(F) = \frac{c_1(\Tvert)\cdot \beta_F+ke}{r}.$$

\noindent
The second cohomology group $H^2(\Grk,{\mathbb Q})$ is generated (over ${\mathbb Q}$) 
by $c_1(\Tvert)$ and the class of a fibre of $\pi$. The class of any curve 
$\beta$ in $H_2(\Grk,{\mathbb Q})$ is therefore determined by the intersection of
$\beta$ with these two classes.  Any section
of $\pi$ will have intersection number one with the class of a fibre, and 
we have seen above that the degree of the quotient bundle $F$ determines the
intersection with $c_1(\Tvert)$ for the corresponding section class $\beta_F$.  
The conclusion is that if we fix the degree of the quotient bundle $d$, then 
there is a {\em unique} class $\beta_d$ such that $\beta_d=\sigma_{F*}[C]$ for
all sections $\sigma_F$ corresponding to quotient bundles $F$ of degree $d$.

\medskip
\noindent
 From the point of view of sections of the relative Grassmannian, one natural way to
compactify the space of quotient bundles of degree $d$ is to look at the Kontsevich
space of stable maps $\Mg(\Grk,\beta_d)$.

\medskip
\noindent
A point of $\Mg(\Grk,\beta_d)$ corresponds to the data of an
isomorphism class of a stable map $(C',f)$, where
$C'$ is a complete connected curve of arithmetic genus $g$, with at most nodes as singularities, 
and $f:C'\rightarrow \Grk$ is 
a map with $f_*[C']=\beta_d$.  The map is {\em stable} in the following sense (see e.g. 
\cite{Fulton} 1.1):
\newline 
-- if a rational component of $C^{\prime}$ is contracted by $f$, then it must contain at least 
three points of intersection with other components
\newline
-- if a component of arithmetic genus $1$ is contracted by $f$, then it must contain at least
one point of intersection with other components (this will not occur in our situation).

\medskip
\noindent
For more details about spaces of stable maps, see 
\cite{Fulton} or \cite{Harris} 2.E, 
and for a construction in arbitrary characteristic see \cite{dan}. 
In the case of characteristic $p>0$, the homology group
$H_2(\Grk,\QQ)$ should be replaced by the 
\'etale cohomology group $H^{2m}_{\mbox{\tiny\'et}}(\Grk,\QQ_{\ell})$ for some prime
$\ell\neq p$ (with $m=k(r-k)=\dim\Grk-1$), and $\beta_d$ by the corresponding cycle class
in that cohomology group.
 
\medskip
\noindent
The fact that $\beta_d$ is the class of a section, and that $C'$ is of the same
arithmetic genus as $C$ puts strong conditions on $C'$ and on the map $f$.  First of
all there must be one component, $C'_0$,  of $C'$ mapping isomorphically to $C$ under the
composition $\pi\circ f$.   The condition on the genus now means that the other
components of $C'$ can only have genus $0$, forming trees of rational curves hanging off of this 
distinguished component.
The map $f$ must send $C'_0$ to a 
section over $C$, and map the rational tails into the Grassmannians $\Gr(E_{p_i},k)$ over
various points $p_i$ of $C$.  Also, \emph{it cannot completely collapse any tree of rational curves} (by the
stability condition) but it can map a rational component multiply onto its image.
A basic picture of one of these maps is shown below.

$${
\setlength{\unitlength}{0.004ex}
\begingroup\makeatletter\ifx\SetFigFont\undefined
\def\x#1#2#3#4#5#6#7\relax{\def\x{#1#2#3#4#5#6}}%
\expandafter\x\fmtname xxxxxx\relax \def\y{splain}%
\ifx\x\y   
\gdef\SetFigFont#1#2#3{%
  \ifnum #1<17\tiny\else \ifnum #1<20\small\else
  \ifnum #1<24\normalsize\else \ifnum #1<29\large\else
  \ifnum #1<34\Large\else \ifnum #1<41\LARGE\else
     \huge\fi\fi\fi\fi\fi\fi
  \csname #3\endcsname}%
\else
\gdef\SetFigFont#1#2#3{\begingroup
  \count@#1\relax \ifnum 25<\count@\count@25\fi
  \def\x{\endgroup\@setsize\SetFigFont{#2pt}}%
  \expandafter\x
    \csname \romannumeral\the\count@ pt\expandafter\endcsname
    \csname @\romannumeral\the\count@ pt\endcsname
  \csname #3\endcsname}%
\fi
\fi\endgroup
{\renewcommand{\dashlinestretch}{30}
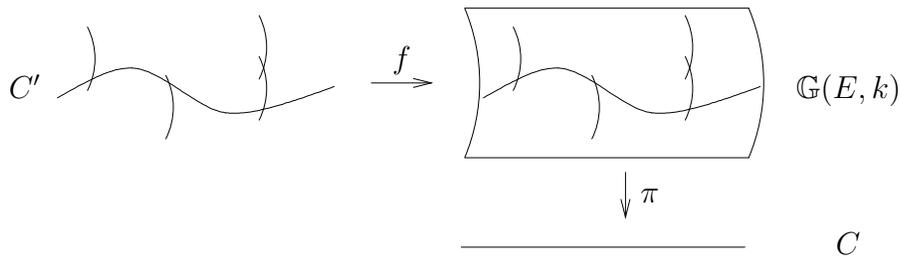
\begin{picture}(15499,4881)(0,-10)
\put(-300.000,3820.500){\arc{2850.987}{5.8195}{6.7468}}
\put(3150.000,3220.500){\arc{2850.987}{5.8195}{6.7468}}
\put(10650.000,3333.000){\arc{7800.000}{5.8884}{6.6780}}
\put(4950.000,3333.000){\arc{7800.000}{5.8884}{6.6780}}
\put(8250.000,3820.500){\arc{2850.987}{5.8195}{6.7468}}
\put(11700.000,3220.500){\arc{2850.987}{5.8195}{6.7468}}
\put(1275.000,2845.500){\arc{2850.987}{5.8195}{6.7468}}
\put(9825.000,2845.500){\arc{2850.987}{5.8195}{6.7468}}
\put(11700.000,4045.500){\arc{2850.987}{5.8195}{6.7468}}
\put(3150.000,4045.500){\arc{2850.987}{5.8195}{6.7468}}
\path(8475,33)(14175,33)
\path(8550,1833)(14250,1833)
\path(8550,4833)(14250,4833)
\path(6675,3333)(7875,3333)
\path(7515.000,3243.000)(7875.000,3333.000)(7515.000,3423.000)
\path(11775,1533)(11775,633)
\path(11685.000,993.000)(11775.000,633.000)(11865.000,993.000)
\path(375,3033)	(448.117,3075.324)
	(518.854,3115.852)
	(587.281,3154.613)
	(653.464,3191.633)
	(717.475,3226.939)
	(779.380,3260.560)
	(839.249,3292.523)
	(897.150,3322.854)
	(953.152,3351.583)
	(1007.324,3378.736)
	(1059.734,3404.340)
	(1110.451,3428.423)
	(1159.544,3451.013)
	(1207.081,3472.137)
	(1253.131,3491.823)
	(1297.762,3510.097)
	(1341.044,3526.988)
	(1383.045,3542.523)
	(1423.834,3556.729)
	(1463.479,3569.634)
	(1502.048,3581.266)
	(1539.612,3591.651)
	(1611.994,3608.792)
	(1681.174,3621.278)
	(1747.703,3629.328)
	(1812.128,3633.162)
	(1875.000,3633.000)

\path(1875,3633)	(1929.016,3628.899)
	(1984.111,3620.369)
	(2040.207,3607.698)
	(2097.226,3591.173)
	(2155.089,3571.084)
	(2213.718,3547.716)
	(2273.034,3521.360)
	(2332.960,3492.302)
	(2393.416,3460.830)
	(2454.324,3427.233)
	(2515.606,3391.797)
	(2577.184,3354.812)
	(2638.978,3316.565)
	(2700.911,3277.343)
	(2762.904,3237.436)
	(2824.879,3197.130)
	(2886.757,3156.714)
	(2948.460,3116.475)
	(3009.909,3076.702)
	(3071.026,3037.682)
	(3131.734,2999.703)
	(3191.952,2963.053)
	(3251.603,2928.021)
	(3310.608,2894.893)
	(3368.890,2863.958)
	(3426.369,2835.504)
	(3482.967,2809.819)
	(3538.606,2787.191)
	(3593.207,2767.906)
	(3646.692,2752.255)
	(3698.982,2740.523)
	(3750.000,2733.000)

\path(3750,2733)	(3790.969,2729.227)
	(3832.685,2726.516)
	(3875.246,2724.890)
	(3918.753,2724.374)
	(3963.305,2724.992)
	(4009.001,2726.767)
	(4055.941,2729.724)
	(4104.225,2733.886)
	(4153.953,2739.279)
	(4205.223,2745.925)
	(4258.135,2753.849)
	(4312.789,2763.075)
	(4369.285,2773.627)
	(4427.722,2785.529)
	(4488.199,2798.805)
	(4550.816,2813.479)
	(4615.673,2829.575)
	(4682.870,2847.117)
	(4752.505,2866.129)
	(4824.679,2886.635)
	(4899.490,2908.659)
	(4937.917,2920.248)
	(4977.040,2932.226)
	(5016.872,2944.595)
	(5057.426,2957.359)
	(5098.714,2970.520)
	(5140.749,2984.082)
	(5183.543,2998.047)
	(5227.108,3012.419)
	(5271.457,3027.201)
	(5316.603,3042.395)
	(5362.558,3058.005)
	(5409.333,3074.033)
	(5456.943,3090.483)
	(5505.399,3107.357)
	(5554.713,3124.660)
	(5604.898,3142.392)
	(5655.967,3160.559)
	(5707.932,3179.162)
	(5760.805,3198.204)
	(5814.599,3217.690)
	(5869.327,3237.620)
	(5925.000,3258.000)

\path(8925,3033)	(8998.117,3075.324)
	(9068.854,3115.852)
	(9137.281,3154.613)
	(9203.464,3191.633)
	(9267.475,3226.939)
	(9329.380,3260.560)
	(9389.249,3292.523)
	(9447.150,3322.854)
	(9503.152,3351.583)
	(9557.324,3378.736)
	(9609.734,3404.340)
	(9660.451,3428.423)
	(9709.544,3451.013)
	(9757.081,3472.137)
	(9803.131,3491.823)
	(9847.763,3510.097)
	(9891.044,3526.988)
	(9933.045,3542.523)
	(9973.834,3556.729)
	(10013.479,3569.634)
	(10052.048,3581.266)
	(10089.612,3591.651)
	(10161.994,3608.792)
	(10231.174,3621.278)
	(10297.703,3629.328)
	(10362.128,3633.162)
	(10425.000,3633.000)

\path(10425,3633)	(10479.016,3628.899)
	(10534.111,3620.369)
	(10590.207,3607.698)
	(10647.226,3591.173)
	(10705.089,3571.084)
	(10763.718,3547.716)
	(10823.034,3521.360)
	(10882.960,3492.302)
	(10943.416,3460.830)
	(11004.324,3427.233)
	(11065.606,3391.797)
	(11127.184,3354.812)
	(11188.978,3316.565)
	(11250.911,3277.343)
	(11312.904,3237.436)
	(11374.879,3197.130)
	(11436.757,3156.714)
	(11498.460,3116.475)
	(11559.909,3076.702)
	(11621.026,3037.682)
	(11681.734,2999.703)
	(11741.952,2963.053)
	(11801.603,2928.021)
	(11860.608,2894.893)
	(11918.890,2863.958)
	(11976.369,2835.504)
	(12032.967,2809.819)
	(12088.606,2787.191)
	(12143.207,2767.906)
	(12196.692,2752.255)
	(12248.982,2740.523)
	(12300.000,2733.000)

\path(12300,2733)	(12340.969,2729.227)
	(12382.685,2726.516)
	(12425.246,2724.890)
	(12468.753,2724.374)
	(12513.305,2724.992)
	(12559.001,2726.767)
	(12605.941,2729.724)
	(12654.225,2733.886)
	(12703.953,2739.279)
	(12755.223,2745.925)
	(12808.135,2753.849)
	(12862.789,2763.075)
	(12919.285,2773.627)
	(12977.722,2785.529)
	(13038.199,2798.805)
	(13100.816,2813.479)
	(13165.673,2829.575)
	(13232.870,2847.117)
	(13302.505,2866.129)
	(13374.679,2886.635)
	(13449.490,2908.659)
	(13487.917,2920.248)
	(13527.040,2932.226)
	(13566.872,2944.595)
	(13607.426,2957.359)
	(13648.714,2970.520)
	(13690.749,2984.082)
	(13733.543,2998.047)
	(13777.108,3012.419)
	(13821.457,3027.201)
	(13866.603,3042.395)
	(13912.558,3058.005)
	(13959.333,3074.033)
	(14006.943,3090.483)
	(14055.399,3107.357)
	(14104.713,3124.660)
	(14154.898,3142.392)
	(14205.967,3160.559)
	(14257.932,3179.162)
	(14310.805,3198.204)
	(14364.599,3217.690)
	(14419.327,3237.620)
	(14475.000,3258.000)

\put(12075,933){\makebox(0,0)[lb]{\smash{{{\SetFigFont{12}{14.4}{rm}$\pi$}}}}}
\put(16000,-100){\makebox(0,0)[lb]{\smash{{{\SetFigFont{12}{14.4}{rm}$C$}}}}}
\put(7100,3633){\makebox(0,0)[lb]{\smash{{{\SetFigFont{12}{14.4}{rm}$f$}}}}}
\put(-600,3033){\makebox(0,0)[lb]{\smash{{{\SetFigFont{12}{14.4}{rm}$C'$}}}}}
\put(15200,3033){\makebox(0,0)[lb]{\smash{{{\SetFigFont{12}{14.4}{rm}$\Grk$}}}}}
\end{picture}
}
}$$
The previous description is slightly more confusing than it needs to be.  
Nodes connecting two rational curves
can always be smoothed, owing to the fact that the tangent bundle to the 
Grassmannian is generated by global 
sections. This means that the general point of any component of 
$\Mg(\Grk,\beta_d)$ consists of a section $C_0'$
mapping isomorphically to $C$, along with a certain number of rational tails, 
each one hanging off of $C_0'$.  None
of these tails are collapsed under the map to the relative Grassmannian, 
and if the fibres of the relative Grassmannian
have dimension greater than two the map will be an immersion at the general point of the component.

\medskip
\noindent
On each $\Gr(E_{p_i},k)$ above a point $p_i$ of $C$, the vertical tangent bundle $\Tvert$ 
restricts to the tangent bundle of this Grassmannian.  The determinant line bundle of this tangent
bundle on the fibre is {\em ample} (more precisely, it is $r$ times the line bundle corresponding 
to the Pl\"ucker embedding).
If $C'_j$ is any rational component of $C'$, this shows that $f_{*}[C'_j]\cdot c_1(\Tvert)$ is
a {\em positive} multiple of $r$.  Since $f_{*}[C']\cdot c_1(\Tvert)=(rd-ke)$, this
means that 
$$f_{*}[C'_0]\cdot c_1(\Tvert) \le (r(d-1)-ke),$$ 
and so the quotient bundle $F_0$ represented by
the section $C'_0$ is of degree {\em less} than $d$.

\medskip
\noindent
As a conclusion, from the discussion above we get the following useful

\medskip
\noindent
{\bf Observation:}
The boundary points of 
$\Mg(\Grk, \beta_d)$ are constructed from 
sections $\sigma_{F_0}$ corresponding 
to quotient bundles of degree {\em strictly smaller} than $d$, along with attached rational
tails which function to increase the total
degree of intersection with $c_1(\Tvert)$.

\medskip
\noindent
{\bf Remark.}
The terminology that a stable map with reducible domain is a {\em boundary point} is standard but potentially slightly
misleading.  It is possible (and common) that there are components of $\Mg(\Grk,\beta_d)$ such that 
every point of the component is a boundary point.  We will also use the name {\em reducible map} for 
such a point.

\medskip
\noindent
For later reference, we end this section by recalling the existence of the more general moduli space
$\Mgn(\Grk, \beta_{d})$ of stable maps with $n$ marked points (see \cite{Fulton} \S1). This comes with a 
natural forgetful morphism to $\Mg(\Grk, \beta_d)$.

\bigskip

\bigskip
\noindent
\mysection{Dimension Estimates for Spaces of Stable Maps}
\label{map_bound}

\medskip
\noindent
For each fixed rank $k$, there is a lower bound on the degrees $d$ of possible quotient bundles $F$ of that rank.  For instance, it
is a simple consequence of Riemann-Roch that $d\ge (e+(r-k)(1-g)-h^0(C,E))$; if $F$ were a quotient bundle of smaller degree, then the
kernel $S_F$ of the map to $F$ would have at least a 
$$ h^0(C,S_F) \ge (e-d) + (r-k)(g-1) > h^0(C,E)$$
dimensional space of sections, which is a contradiction.

\medskip
\noindent
For each rank $k$, let 
$$d_k = d_k (E):=\underset{{\rm rk}F=k}{{\rm min}}\left\{{{\rm deg}F~\rule{0cm}{0.6cm}|\,\,\,\mbox{$F$ a quotient of $E$}}\right\}.$$  
This number will depend on the vector bundle $E$, and is lower semi-continuous
in families of vector bundles.

\medskip
\noindent
The purpose of this section is to prove the dimension estimate:

\begin{theorem}\label{4}\label{Mg_dim}
$$ \dim\Mg(\Grk,\beta_d) \le k(r-k) + (d-d_k)r, \mbox{ for all $d\ge d_k$.}$$
\end{theorem}

\medskip
\noindent
This estimate includes those components of $\Mg(\Grk,\beta_d)$ whose general points do not correspond to quotient
bundles. For these the inequality is actually strict,
as the proof will show. 

\bigskip
\noindent
The method of obtaining the bound is very simple.  We will show 
how the \emph{bend and break} lemma of Mori guarantees 
that once any component of 
$\Mg(\Grk,\beta_d)$ has dimension strictly greater than $k(r-k)$, there
is a divisor in that component corresponding to \emph{reducible maps} (i.e., 
stable maps with reducible domain curves).  From the
previous analysis, these correspond to quotients of lower degree, and this
will inductively give us the dimension bound.  This application of bend and
break is the content of the following lemma.

\medskip
\noindent
\begin{lemma}\label{5}
\label{breaking_lemma}
If $X$ is a component of $\Mg(\Grk,\beta_d)$ whose generic point corresponds to a map with irreducible domain, 
and $\dim\,X > k(r-k)$, then there is a (nonempty) divisor $Y\subset X$ consisting of reducible maps.
\end{lemma}

\medskip
\noindent
Before proving the lemma, let us use it to prove the dimension bound.

\medskip
\noindent
{\em Proof of Theorem \ref{4}} (using lemma):
First of all, suppose that $d=d_k$.  Every point of $\Mg(\Grk,\beta_{d_k})$ corresponds to an irreducible map -- 
since a reducible map
would give a quotient bundle of strictly smaller degree.  The lemma then shows that no component
of $\Mg(\Grk,\beta_{d_k})$ can have dimension greater than $k(r-k)$ (the key is that $Y$ is nonempty).  This gives
the base case for the induction.

\medskip
\noindent
For $d>d_k$, suppose that $X$ is a component of $\Mg(\Grk,\beta_d)$ satisfying the conditions of the lemma.  
The divisor $Y$ produced by the lemma
consists of reducible maps.  Let $Y'$ be a component of $Y$ of the same dimension.  
By the discussion in the previous section,
we can assume that the generic point of $Y'$  consists of a
reducible map 
where the component $C_0'$ corresponds to a quotient bundle of degree 
$(d-\delta)$, along with a rational component
 which makes up
the difference when intersected with $c_1(T_V)$; i.e.,  if $C_1'$ is the rational tail, 
then $f_{*}[C_1']\cdot c_1(T_V) = r\delta$.

\medskip
\noindent
This set of reducible maps can be built up in the usual way by using smaller spaces of stable maps, and this will give us the
inductive dimension bound.

\medskip\noindent
Let $\alpha_\delta$ be the homology class in $H_2(\Grk, \ZZ)$ corresponding to curves of the same type as these rational tails -- these are curves
contained in a fibre of the map $\pi:\Grk\longrightarrow C$, and which have intersection $r\delta$ with the vertical tangent bundle.
 The space of maps $\Moi(\Grk,\alpha_\delta)$ of 
rational tails with a marked point is of dimension $r\delta + k(r-k) -1$.

\medskip
\noindent
The component $Y'$ of $Y$ is in the image of
$$\Mgi(\Grk,\beta_{d-\delta}) \times_{\Grk} \Moi(\Grk,\alpha_\delta) \longrightarrow \Mg(\Grk,\beta_d)$$
where we glue together the maps whose marked points are sent to the same point of $\Grk$.  This fibre product 
has dimension
\begin{equation*}
\begin{split}
& \dim\,\Mg(\Grk,\beta_{d-\delta}) + 1 + \dim\,\Moi(\Grk,\alpha_\delta) -\dim\,\Grk \\
& = \dim\,\Mg(\Grk,\beta_{d-\delta}) + 1 + r\delta + k(r-k) -1 - k(r-k) -1 \\
& = \dim\,\Mg(\Grk,\beta_{d-\delta}) + r\delta -1.
\end{split}
\end{equation*}

\noindent
Since $Y'$ is in the image of this map, the dimension of $Y'$ is less than or equal to the dimension above.  
Hence:
$$\dim\,X = \dim\,Y' + 1 \le \dim\,\Mg(\Grk,\beta_{d-\delta}) + r\delta. $$
Using the inductive dimension bound for the dimension of $\Mg(\Grk,\beta_{d-\delta})$ this gives us
$$\dim X \le k(r-k) + (d-\delta-d_k)r + r\delta = k(r-k) + (d-d_k)r,$$
which is exactly what we wanted.

\medskip
\noindent
This shows why the lemma gives us the dimension bound, at least as long as we consider components $X$ of $\Mg(\Grk,\beta_d)$ 
whose generic map is irreducible.  In the case that $X$ is a component whose generic point {\em does} correspond to a reducible map, 
we can express $X$ (as above) as being in the image of a product of spaces parametrizing rational tails, and a space $\Mg(\Grk,\beta_{d'})$ 
where $d'$ is of smaller degree.  If the generic point of $X$ consists of a map with $l$ rational tails, then the same kind of 
argument as above shows that
$$\dim\,X \le k(r-k) + (d-d_k)r - l.$$
This finishes the proof of the theorem using 
Lemma \ref{5} \doneproof.

\bigskip
\noindent
The only step remaining is to prove the lemma.

\medskip
\noindent
{\em Proof of Lemma:} Let $X$ be a component satisfying the conditions of the lemma.  It is easier to pass to the space $\Mgi(\Grk,\beta_d)$
with a single marked point;  let $\widetilde{X}$ be the inverse image of $X$ under the forgetful map $\Mgi(\Grk,\beta_d)\longrightarrow
\Mg(\Grk,\beta_d)$.  It is sufficient to prove that there is a nonempty divisor $\widetilde{Y}$ in $\widetilde{X}$ consisting of 
reducible maps.  (In fact it would be sufficient to prove that this component has at least one point corresponding
to a reducible map, since the locus of nodal curves in a family is always of codimension at most one, provided it
is nonempty.)

\medskip
\noindent
The advantage of passing to the space with a marked point is that we have an evaluation map $ev:\widetilde{X} \longrightarrow \Grk$.
By hypothesis, $\widetilde{X}$ has dimension strictly greater than $k(r-k) + 1$.  Since $\Grk$ is of dimension $k(r-k)+1$, we see
that all the fibres of the evaluation map are at least one dimensional.  To show the existence (and dimension) of $\widetilde{Y}$, it
is sufficient to show that each fibre has a nonempty set of reducible maps, of codimension zero or one.

\medskip
\noindent
Pick a point $p$ of $C$, and let $F_p$ be a $k$ dimensional quotient of the vector space $E_p$.  Let $Z_{F_p}=ev^{-1}([F_p])$ be the fibre of
the evaluation map over this point in $\Grk$.  If $Z_{F_p}$ is empty there is nothing to prove.  If $Z_{F_p}$ consists completely of 
reducible maps,
then there is still nothing to prove.  The remaining possibility, which has to happen generically over the image, is that $Z_{F_p}$ has some points which correspond to maps with an 
irreducible
domain curve.  If the locus of reducible maps is {\em not} a (nonempty) divisor, then we can find an irreducible complete curve $B''$ in 
$Z_{F_p}$ such that every point of $B''$ corresponds to a map with 
irreducible domain curve.  (Also note, by the discussion in \S2, 
that these
curves are all isomorphic to $C$).  It is now easy to show that the existence of $B''$ leads to a contradiction.

\medskip
\noindent
Let $B'$ be the normalization of $B''$, so that we can work with a 
smooth base curve.  
We have a map $B' \longrightarrow B'' \longrightarrow \Mgi(\Grk,\beta_d)$.  Since the 
moduli space of maps is only a coarse moduli space, this might not correspond to an actual family.  However, there is a finite cover
$B\longrightarrow B'$ so that the induced map $B\longrightarrow \Mgi(\Grk,\beta_d)$ {\em does} come from a family ${\cal C}$ over $B$.  
This family is isotrivial, since all fibres are isomorphic to $C$. 

\medskip
\noindent
At this stage we have the family ${\cal C}$ over $B$, along with a section $s:B\longrightarrow {\cal C}$ corresponding
to the marked point, and a map $f:{\cal C} \longrightarrow \Grk$ inducing the stable maps on the fibres.  The combination of the map
$f$ followed by the projection  $\pi$ from $\Grk$ to $C$ gives an isomorphism of each fibre of the family with $C$, and 
so we see that ${\cal C}$ is the trivial family ${\cal C}=C\times B$ over $B$.  

\medskip
\noindent
The section $s$ is sent to the point $F_p$ under $f$, since $B''$ is
contained in the fibre $Z_{F_p}$ of the evaluation map.   
This shows that the images of the maps $C\rightarrow \Grk$ 
parametrized by $B$ go through the fixed 
point $F_p$. Now this family 
of stable maps over $B$ induced by $f$ is not the constant family, since $B''$ was a nontrivial curve 
in the moduli space.  
This contradicts the bend and break lemma of Koll\'ar and Mori (see \cite{Mori} 1.7), which asserts 
that in the given situation at least one of the image curves should ``break'' and have a 
rational component passing through $F_p$.
We conclude that the curve $B''$ could not exist, and therefore that the locus of reducible 
curves in $Z_{F_{p}}$ is a nonempty divisor \doneproof.

\medskip
\noindent
{\bf Remark.} What we use is actually only the rigidity part of the bend and break lemma,
and not the full strength of Mori's arguments.  The referee has pointed out that the
rigidity statement can be proved directly in this case, as follows:

\medskip
\noindent
The family $C\times B \longrightarrow \Grk$ gives a short exact sequence of 
bundles over $C\times B$

$$ 0\longrightarrow S \longrightarrow \pr_1^{*}(E)\longrightarrow Q \longrightarrow 0.$$

\medskip
\noindent
For the special point $p$ of $C$, the hypothesis is that $Q$ restricted to
$\{p\}\times B$ is the trivial bundle, which implies that $Q$ has degree
zero when restricted to $\{q\}\times B$ for any point $q\in C$.

\medskip
\noindent
For any such point $q$, the restriction of $\pr_1^{*}(E)$ to $\{q\}\times B$ is 
the trivial bundle.  The only degree zero quotients of the trivial bundle are themselves trivial
bundles, {\em via trivial quotients}, i.e., quotients which do not depend on the point of $B$.

\medskip
\noindent
This means that for each $q\in C$, the family of quotients parameterized by $B$ is constant,
and therefore that the whole family of quotients on $C$ does not vary over $B$, contrary to the
way that $B$ was initially chosen.

\medskip
\noindent
{\bf Remark.} Often in arguments using degeneration the difficult part is 
proving that a one parameter family actually has to degenerate.
In our situation, the degeneration argument is very easy,
the advantage coming from the fact that the curves in our families do not
vary in moduli.

\pagebreak
\bigskip
\noindent
{\bf Some variations on this bound.}
Let $m_k$ be the dimension of the space of stable maps in the case of rank $k$ 
quotients of minimal degree.  i.e.:
$$m_k := {\rm dim~} \Mg(\Grk, \beta_{d_k}).$$ 
In order to use induction in the proof of Theorem \ref{4}, we need that at
each stage the dimension of the component we are considering be greater than
$k(r-k)$.
Therefore, whenever $m_k \ge k(r-k) - r$, the same argument also gives the 
slightly better result
\begin{equation}
\dim \Mg (\Grk, \beta_d)\leq m_k + (d-d_k)r.
\end{equation} 

\medskip
\noindent
There is also a lower bound on the dimension of $\Mg(\Grk,\beta_d)$ 
involving the same number $m_k$ 
(and which does not require that $m_k\ge k(r-k)-r$).  
The bound is obtained by taking sections $C_0'$ 
corresponding to minimal quotients and gluing on a single rational curve which
maps to one of the fibres, with the degree of the rational curve chosen to 
make up the difference between $d$ and $d_k$.  There are $m_k$ parameters
for the choice of $C_0'$, $1$ parameter for the choice of where to glue the
rational curve onto $C_0'$, and $(d-d_k)r-2$ parameters for the choice of the
map of the rational curve into the Grassmannian fibre.  This gives:

\begin{proposition}\label{20}
${\rm dim~} \Mg (\Grk, \beta_d)\geq m_k + (d-d_k)r - 1$.
\end{proposition}

\bigskip
\noindent
{\bf Remark.} Going back to the case $m_k \ge k(r-k) - r$, the upper bound there
is one more then this lower bound, and this leaves only two 
possibilities for the dimension of $\Mg(\Grk,\beta_d)$.  Which possibility 
occurs depends on whether or not it is possible to smooth the nodes in the 
families above.

\medskip
\noindent
The same trick of gluing on a rational curve to make up the difference in
degree gives an
inequality between the dimensions of any two of the stable map spaces,
namely
$$ \dim \Mg(\Grk,\beta_{d_2}) \ge \dim \Mg(\Grk, \beta_{d_1}) + (d_2-d_1)r -1,\,\,
\mbox{for $d_2\ge d_1$}.$$

\medskip
\noindent
Combining this information about the growth of the dimensions with the same
inductive argument of Theorem \ref{4} gives the following other variation:
If $d'$ is the first degree so that dim $\Mg(\Grk,\beta_{d'}) \geq k(r-k)-r$, then

\begin{equation}
\dim \Mg(\Grk,\beta_d) \le k(r-k) + (d-d')r, \,\mbox{for all $d\ge d^{\prime}-1$}.
\end{equation}

\medskip
\noindent
Our basic picture in these upper (and lower) bounds is that the dimension
grows by $r$ when the degree increases by $1$, at least once the dimension
of the space is large enough.  It is tempting to wonder if this is always
the case (even when the dimensions are small), and if, for example, inequality
$(1)$ might be true in all cases.  This is in fact not true, as Example 
5.\ref{formula_not_always_true} will show. 

\bigskip

\pagebreak
\bigskip
\noindent
\mysection{Dimension Estimates for Quot Schemes}

\medskip
\noindent
We now switch to looking at the other, more familiar, compactification
of the space of vector bundle quotients of a fixed vector bundle $E$, namely 
Grothendieck's Quot scheme.
Let ${\rm Quot}_{k,d}(E)$ be the the scheme parametrizing coherent quotients
of $E$ of generic rank $k$ and degree $d$.  
A general introduction to the theory of Quot schemes is given e.g. in \cite{Le Potier} \S4 and \S8. 
We would like to prove the identical dimension estimate:

\begin{theorem}\label{1}\label{Quot_dim}
$$\dim {\rm Quot}_{k,d}(E) \leq k(r-k)+ (d-d_{k})r, \mbox{ for all $d\ge d_k$.}$$
\end{theorem}

\medskip
\noindent
{\em Proof.}
If $Q$ is a component of $\Quot$ whose generic point parametrizes vector
bundle quotients, then there is nothing to prove, since it is birational
to a component of $\Mg(\Grk,\beta_d)$ and Theorem \ref{4}
applies.  As in the case of $\Mg(\Grk,\beta_d)$, however, there are usually
components whose general point does not correspond to a locally free quotient,
and, as in the case of $\Mg(\Grk,\beta_d)$, we will show that 
for these the inequality is always strict.

\medskip
\noindent
Let $Q$ be a component of this type. Every non locally free quotient 
$F$
determines a diagram:
$$\xymatrix{
& & 0 \ar[d] & 0 \ar[d] \\
0 \ar[r] & S_F \ar[r] \ar[d]_{\cong} & S_{F^\prime} \ar[r] \ar[d] & \tau \ar[r]
\ar[d] & 0 \\
0 \ar[r] & S_F \ar[r] & E \ar[r] \ar[d] & F \ar[r] \ar[d] & 0 \\
& & F^{\prime} \ar[d] \ar[r]^{\cong} & F^{\prime} \ar[d] \\
& & 0 & 0 } $$
where $S_{F^\prime}$ is the saturation of the kernel $S_F$, 
$F^{\prime}$ is a quotient vector bundle of degree $(d-a)$ and $\tau$, the torsion
subsheaf of $F$, is given by a 
nontrivial zero-dimensional subscheme, say of length $a>0$. We can
stratify the set of all such $F$'s in $Q$ 
according to the value of the parameter $a$, which runs over a finite set. If we denote
by $\{F\}_{a}$ the subset of $Q$ corresponding to a fixed $a$, this gives:
$${\rm dim~}\{F\}_{a}\leq {\rm dim~}{\rm Quot}_{k,d-a}(E)+(r-k)a.$$
Assuming the dimension bound true for smaller $d$, we get 
the inequality:
$${\rm dim}\{F\}_{a}\leq k(r-k)+(d-a-d_{k})r+(r-k)a=$$
$$=k(r-k)+(d-d_{k})r-ka.$$
Letting $a$ run over all possible values (or the minimum value) which occur 
for the component $Q$ finishes the proof of the bound \doneproof.

\bigskip
\noindent
{\bf Remark.} With only slight rewording the proof also gives the stronger
$$\dim \Quot \le \dim \Mg(\Grk,\beta_d) \quad\,\mbox{for all $d\ge d_k$}.$$

\medskip
\noindent
As in the case of $\Mg (\Grk, \beta_d)$ (see \ref{20}), we can easily produce a lower 
bound for the dimension of $\Quot$, this time by starting with vector bundle quotients of 
lower degree and going
up via elementary transformations. Again, it may well be the case that the components achieving 
such a bound consist entirely of non-locally free quotients.

\medskip
\noindent
\begin{proposition}\label{21}
Let $d\geq d^\prime \geq d_k$ and $\delta = (d - d^\prime)$. Then we have 
$${\rm dim~} {\rm Quot}_{k,d}(E) \geq {\rm dim~} {\rm Quot}_{k,d^\prime}(E) + \delta (r-k).$$
In particular
${\rm dim~} {\rm Quot}_{k,d}(E) \geq m_k + (d-d_k)(r-k)$, 
 for all $d\geq d_k$.
\end{proposition}

\medskip
\noindent
{\em Proof.}
We only need to consider the diagram in the proof of Theorem \ref{1} when 
$F^{\prime}$ 
runs over the family of quotients of degree $d^\prime$. This gives at most 
${\rm dim~} {\rm Quot}_{k,d^\prime}(E)$ parameters, and
the other $\delta (r-k)$ parameters come from the fact that $S_F$ is an (inverse)
elementary transformation of $S_{F^\prime}$ of length $\delta$. The last statement 
follows by making $d^\prime = d_k$.  

\bigskip
\noindent
{\bf Remark.} 
Propositions \ref{20} and \ref{21} show that in general there is more 
flexibility 
in dealing with stable maps with reducible domain rather than non-locally free 
quotients.   In section \S7 we will show that there is always a surjective
morphism $\Mg(\Grk,\beta_d)\longrightarrow \Quot$ extending the map on
smooth curves if $k=(r-1)$, and no such map if $k\ne(r-1)$.  In general the
relation between the two compactifications seems a bit mysterious.

\bigskip

\bigskip
\noindent
\mysection{Examples and Optimality of the Bound}
\label{examples}

\medskip
\noindent
\example \label{trivial_example} 
Let $E$ be the trivial bundle of rank $r$ on $C$, and let $k$ be any number 
from $1$ to $(r-1)$.  A ``section of the  
relative Grassmann bundle'' in this case is just the graph of a map from 
$C$ to the usual Grassmannian $\Gr(r,k)$. The degree of the corresponding quotient bundle 
is the degree of the linear series associated to the map from $C$ 
to $\PP^{N}\,\left(\mbox{with}\,N={{{r}\choose{k}} -1}\right)$
induced by the Pl\"{u}cker embedding.

\medskip
\noindent
The minimal degree of a quotient bundle is zero (i.e. $d_k=0$ for all $k$) and
in this case we have $\dim \Mg(\Grk,\beta_0)=k(r-k)$.  For $d=1$, there will
be no vector bundle quotients unless the genus of the curve is zero, since
for $g\ge 1$ there are no maps from $C$ to $\PP^N$ of degree $1$.  If $g\ge 1$,
all the points of $\Mg(\Grk,\beta_1)$ correspond to maps mapping $C$ to
a point, and mapping an attached rational curve in as a line.  In the Quot
scheme case, all the degree $1$ quotients are torsion.

\medskip
\noindent
For degree $d=2$, there will again be no vector bundle quotients, unless the
curve is hyperelliptic, but there will always be reducible maps or torsion
quotients as above.  For other low $d$, the geometry of both the stable map
and Quot scheme depend on the particular special linear series that live on
the curve.

\medskip
\noindent
For large $d$, the behaviour of linear series on $C$ is regular, and there will
be components of both spaces corresponding to vector bundle quotients, of
dimension $dr-k(r-k)(g-1)$.  
The behaviour of the largest components, however, are
different in the two cases.
Following the constructions of 
Proposition \ref{20}, $\Mg(\Grk,\beta_d)$ has a locus
of dimension $rd+k(r-k)-1$ constructed by mapping $C$ to a point, and mapping
an attached rational curve as a curve of degree $d$. If $g\ne 0$, this is 
larger than the dimension of points parametrizing vector bundle quotients, and
so constitutes a separate component.
On the other hand, if we apply Proposition \ref{21} 
to the Quot scheme, we get a torsion locus of
dimension $d(r-k) +k(r-k)$.  For large $d$ this is a sublocus of the component
of $\Quot$ whose general point parametrizes vector bundle quotients (see
Theorem \ref{irred_of_Quot}) and not a new component at all.

\medskip
\noindent
This example shows that there are vector bundles $E$ where the bound is exact
for $d=d_k$, and gives an example of a case where the families constructed
in Proposition \ref{20} form components of the moduli spaces,
i.e. where it is impossible to smooth the nodes of the family from Proposition
\ref{20}.

\bigskip
\noindent
\example \label{extension_example}
On a curve $C$ of genus $g\ge 2$, 
let $L$ be a line bundle of degree $1$ with no global section, 
and let $G_1$ be a generic extension of the type
$$0\rightarrow \stackrel{r-1}{\oplus} {\cal O} \rightarrow G_1 \rightarrow 
L \rightarrow 0.$$
If $g\ge (r-1)$ then the usual dimension counting shows that the generic
extension of this type is a stable bundle (see also \cite{BGN}
\S5 for more general counts of this kind).

\medskip
\noindent
The minimum degree of any quotient is $d_k=1$, and for any $k$ from 
$1$ to $(r-1)$ there is at least a $(k-1)(r-k)$ dimensional family of degree
$d=1$ quotients obtained by dividing out by a rank $(r-k)$ subbundle of the 
direct sum (if $g\ge r$ then these are the only quotients of degree $1$).

\medskip
\noindent
In most uses of a bound for the dimension of the quotients of a fixed bundle
$E$, what is really needed is a bound for the {\em moduli} of the 
bundles $F$, i.e. the dimension of the isomorphism classes of $F$, throwing
away the data of the particular morphisms $E\rightarrow F$ which 
express $F$ as a quotient.

\medskip
\noindent
In this example, ($d=1$, $g\ge (r-1)$, and the extension generic) 
it is not hard to check that
all the quotients are nonisomorphic,
and this provides an example of a family of quotients with maximal 
variation in moduli (so that the dimension of the space of quotients is the 
same as the dimension of the moduli of $F$'s) as well as an example of 
a stable bundle where the dimension in the case $d=d_k$ is close 
to the upper bound.

\medskip
\noindent
A similar example is to take a generic extension
$$ 0 \rightarrow L^{*} \rightarrow G_0 \rightarrow G_1 \rightarrow 0$$
with $L^{*}$ the dual of the line bundle $L$ above.  For $g\ge r$ this
extension is also stable, and provides an example of a stable bundle with
trivial determinant with a large family of quotients of degree $1$ 
(which also vary maximally in moduli, since they are the same quotients as
above).

\bigskip
\noindent
\example \label{formula_not_always_true}
For any $r\ge 11$, pick $k$ with $r/2<k<r$ such that $(r-k)(2k-r-1)>r$ 
(this works out to be any $k$ between $(r+5)/2$ and $(r-2)$).
Let $G_1$ be a bundle of rank $k$ of the type from Example 
\arabic{section}.\ref{extension_example}, and
let $H$ be a stable bundle of rank $(r-k)$ and degree $1$.  
Set 
$$E = H \oplus G_1,$$ 
which is a bundle of rank $r$ and degree $2$.

\medskip
\noindent
By the stability of $G_1$ and $H$, $E$ can have no quotients of degree $0$. 
If $F$ is a quotient bundle of rank $k$ and degree $1$, 
then the kernel bundle $S_F$ fits in the following diagram (where the vertical
arrows are injections):

$$\xymatrix{
0 \ar[r] & S_F^{\prime} \ar[r] \ar[d] & S_F \ar[d] \ar[r] 
& S_F^{\prime\prime} \ar[r] \ar[d] & 0 \\
0 \ar[r] & H \ar[r] & E \ar[r]  & G_1 \ar[r]  &0}$$
This shows that the slopes of $S_F^{\prime}$ and $S_F^{\prime\prime}$
are less than $1/(r-k)$ and $1/k$ respectively.  Since $1/k < 1/(r-k)$ and
the slope of $S_F$ {\em is} $1/(r-k)$, the only possibility is $S_F^{\prime}=H$
and $S_F^{\prime\prime}=0$, so that $G_1$ is the only quotient bundle of
rank $k$ and degree $1$.  Or, in the notation of \S3,
we have $d_k(E)=1$, $m_k(E)=0$.

\medskip
\noindent
At the end of \S3 we wondered if the bound $(1)$
$${\rm dim~} \Mg (\Grk, \beta_d)\leq m_k + (d-d_k)r$$
might be true in all cases.  In this example, for $d=2$, this would give the 
bound 
$${\rm dim~} \Mg (\Grk, \beta_2)\leq r.$$

\medskip
\noindent
One way to get quotient bundles of degree $2$, rank $k$ is by taking quotients
of $G_1$ which have degree $1$ and rank $(2k-r)$.  We know the dimension of
this family by the calculation in
Example \arabic{section}.\ref{extension_example}, and this shows that
${\rm dim~} \Mg (\Grk, \beta_2)\ge (2k-r-1)(r-k)$, which we chose at the 
beginning to be strictly larger than $r$. 

\medskip
\noindent
{\bf Remark.} The idea of the proof of the bound is to use degeneration to
relate bundles of higher degree to bundles of lower degree, which allows us
to see that the dimension goes up by $r$ as the degree increases.  The only
way that $(1)$ might not be true in general is if there were families of 
quotients which have no relation to each other at all, and that is what
this example provides.

\bigskip
\noindent
\example \label{generic_stable_bundle}
Let $E$ be a generic stable bundle of degree $e$ and rank $r$ on a curve
$C$ of genus $g\geq 2$. 
The behaviour of the dimensions of the Quot schemes (which also applies to
$\Mg(\Grk,\beta_d)$) in this situation has been worked out 
in \cite{Teixidor}, and a summary is the following:

\noindent
\begin{itemize}
\item The invariant $d_k$ is the smallest value which makes 
the expression \\
$d_kr-ke-k(r-k)(g-1)$ nonnegative.
\item The number $m_k$ is the value of $d_kr-ke-k(r-k)(g-1)$  \\
(i.e. the residue class
of $k^2(g-1)-ke$ mod $r$).
\item $\dim \Mg(\Grk,\beta_d) = \dim \Quot = dr-ke -k(r-k)(g-1) \\
 = m_k + (d-d_k)r$ for all $d\ge d_k$.
\end{itemize}

\medskip
\noindent
The bundles satisfying the first condition are sometimes called 
\emph{Lange generic} and were first studied in \cite{Lange}.
If a bundle satisfies all three conditions above, we will say that it is
\emph{generic in the sense of example 5.4}, or sometimes just \emph{5.4 stable}.

\medskip
\noindent
The expression $dr-ke-k(r-k)(g-1)$ appearing repeatedly above is 
the cohomological lower bound for a quotient bundle $F$ of rank $k$ and 
degree $d$:
$$\dim \Mg(\Grk,\beta_d) \ge \chi(S_F^{*}\otimes F) = dr-ke-k(r-k)(g-1).$$

\medskip
\noindent
This example shows both that there are situations where this lower bound
is achieved, and where formula $(1)$ in \S3 may be valid without the condition
$m_k\ge k(r-k)-r$. Most importantly, in \S6 we will see that 
these are precisely the bundles for which $\Mg(\Grk,\beta_d)$ is 
irreducible for large $d$.

\medskip
\noindent
For special choice of the parameters, the dimension of $\Mg(\Grk,\beta_d)$
in this example exactly matches our upper bound.
If $k=1$ and $r$ divides $(g-e)$, or $k=(r-1)$ and $r$
divides $(g+e)$, then $\dim \Mg(\Grk,\beta_d) = (r-1) + (d-d_k)r$ for all 
$d\ge d_k$, which is the upper bound (similarly for $\Quot$).  
For all other choices of parameters the upper bound will be strictly larger.

\bigskip
\noindent
{\bf How good is this bound?} Viewed as functions of $d$, 
both our upper bound (\ref{Mg_dim} or \ref{Quot_dim})
and the cohomological lower bound mentioned above are 
linear functions  
with the same slope $r$, and therefore this part of the bound is 
certainly sharp.  
There are also examples (5.\ref{trivial_example} and 
5.\ref{generic_stable_bundle}) 
where the values of the ``constant term'' are achieved and so without
extra conditions on the bundle this cannot be improved either.
Further evidence that the constant term is good is
that the difference between upper and lower bounds gives the inequality:
$$ rd_k-ke \le k(r-k)g,$$
which is a theorem of Lange \cite{Lange} and Mukai-Sakai \cite{Mukai} 
generalizing a rank $2$ result of Segre and Nagata \cite{Nagata}
on minimal sections of ruled surfaces.

\medskip
\noindent
One appealing aspect of the upper bound is that it depends on
very little data.
The rank and degree ($k$ and $d$) of the quotients 
and the rank $r$ of the bundle are the basic numerical invariants of the 
problem, and the only thing which really ties the formula to the particular
vector bundle $E$ is the (admittedly sometimes hard to understand) invariant
$d_k(E)$.  In particular, the formula is independent of the genus of the
curve, the characteristic of the field, or any other 
conditions on the bundle $E$ (including its degree). 

\medskip
\noindent
Even with
extra conditions on the bundle $E$, the only possible improvement is 
better (and somewhat limited) control over the behaviour of the constant term.
On the other hand, it is exactly the lack of extra conditions on $E$ that
is one of the strongest features of this bound. This lack of restrictions 
makes it very useful in computations;  Lemma 6.3 in the next 
section is an illustration of this.

\bigskip
\noindent
\mysection{Components of Quot and $\Mg(\Grk,\beta_d)$ for large $d$}

\bigskip
\noindent
Let $E$ be a bundle of rank $r$ and degree $e$ on a curve of 
genus $g>0$. We
start by proving a slightly weaker version of the theorem we are aiming for.

\begin{proposition}\label{24}\label{good_component}
For all large $d$ there is a unique
component of ${\rm Quot}_{k,d}(E)$ and $\Mg(\Grk,\beta_d)$ 
whose generic point corresponds to a stable bundle quotient; this
component is of dimension $rd - ke - k(r-k)(g-1)$.
\end{proposition}
We will need two vanishing lemmas. 

\noindent
\begin{lemma}
\label{good_homs}
There is a number $d_1=d_1(E)$ such that for all stable bundles $F$ of rank $k$,
degree $d\ge d_1$, $\mathcal{H}om(E,F)$ is generated by its global sections and
$H^1(E^{*}\otimes F)=0$.
\end{lemma}

\medskip
\noindent
{\em Proof.} Let $L$ be a line bundle of degree $1$.  For any bundle 
$F^{\prime}$, 
there is a number $n$ so that $E^{*}\otimes F^{\prime} \otimes L^{\otimes m}$ 
is globally generated, and so that $H^1(E^{*}\otimes F^{\prime} 
\otimes L^{\otimes m})$ 
vanishes for all $m\ge n$.

\medskip
\noindent
The condition we want is open, and so given any family $X$ of vector bundles
$F^{\prime}$, the data of the minimum possible $n_x$ for each 
member $F^{\prime}_x$ of the family
stratifies $X$ into locally closed subsets.  In particular, this function
has a maximum value on $X$, and so taking $n$ equal to this shows we can
choose an $n$ which works for all members of the family.

\medskip
\noindent
Apply this to the moduli spaces $U_C^s(k,j)$ of rank $k$ stable 
bundles for degrees $j=0,\ldots,(k-1)$, and choose an $n$ which works for
all of these families. Set $d_1=nk$.

\medskip
\noindent
Any stable bundle $F$ of rank $k$ and degree $d\ge d_1$ can be written
as $F=F^{\prime}\otimes L^{\otimes m}$ for $m\ge n$ (the $n$ above) with 
$F'$ a stable bundle of rank $k$ and degree between $0$ and $(k-1)$, 
which proves the lemma \doneproof.

\medskip
\noindent
{\bf Remark.} Considering filtrations shows that it is
possible to choose a $d_1$ so that Lemma \ref{good_homs} holds
for all {\em semi}-stable bundles $F$ of rank $k$, degree $d\ge d_1$.

\noindent
\begin{lemma}
\label{smooth_point}
There is a number $d_2=d_2(E)$ so that for any quotient bundle $F$ of rank
$k$, degree $d\ge d_2$ with either $F$ or $S_F$ $($the kernel$\,)$ semistable,
 $h^{1}(S_F^{*}\otimes F) = 0$
$\,($and so $h^{0}(S_F^{*}\otimes F) = rd - ke - k(r-k)(g-1))$. 
\end{lemma}

\medskip
\noindent
{\em Proof.}
Since the proofs in the two cases are very similar, we only look at
the case when $F$ 
(as opposed to $S_F$) is semistable.
In this case, the result follows immediately from Lemma \ref{good_homs}
and the long exact sequence of cohomology, if we choose $d_2$ to be 
equal to the $d_1$ we obtain in the above Remark
\doneproof.

\bigskip
\noindent
{\em Proof of Proposition \ref{good_component}.}
Let $F$ be a stable bundle of rank $k$ and 
degree $d\ge\max(d_1(E),d_2(E))$. Lemma \ref{good_homs} 
guarantees
that $\mathcal{H}om(E,F)$ is globally generated, and since $k<r$ it is not hard to
see that the generic homomorphism is a surjection, which means that
both $\Quot$ and $\Mg(\Grk,\beta_d)$ contain a point corresponding to the
quotient
$$ 0\longrightarrow S_F \longrightarrow E \longrightarrow F \longrightarrow 0.$$
Lemma \ref{smooth_point} now says that this point is a smooth
point of the component it lies on, and that the dimension of this component
is $rd-ke-k(r-k)(g-1)$.  Since a small 
deformation of a stable bundle is again a stable bundle, the general point of 
this component also corresponds to a stable quotient, and this proves all the
statements of the proposition except for the uniqueness. 

\medskip
\noindent
By Lemma \ref{good_homs} $h^1(E^{*}\otimes F)=0$ for all $F$ in the moduli
space $U_C^s(k,d)$ of stable rank $k$ degree $d\ge d_1$ bundles, and therefore
the vector spaces $\Hom(E,F)$ form a bundle over $U_C^s(k,d)$.  Since $U_C^s(k,d)$
is irreducible, the total space of the bundle is as well. (This should
be interpreted correctly in the usual way: even though there might not be a
Poincar\'e bundle over the moduli space, things are fine at least after passing 
to an \'etale cover.) The points 
corresponding to surjections from $E$ onto a stable bundle $F$ of rank $k$,
degree $d$ form an open subset of this total space.  The points of $\Mg(\Grk,
\beta_d)$ and $\Quot$ corresponding to stable quotients are in one to one correspondence
with the points of this irreducible open set, and this proves uniqueness
\doneproof.

\medskip
\noindent
{\bf Remark.} The proof shows in particular that, as long as 
$d$ is large enough, \emph{every} semistable bundle 
of degree $d$ is a quotient of $E$ and moreover, all such 
quotients correspond to points in the smooth locus.

\bigskip
\noindent
The uniqueness in Proposition \ref{good_component}
is true even if we 
ask just that the general point of the component parametrize vector bundle
quotients (without necessarily being stable).
To prove this we will use the proposition
and the following lemma on avoiding quotients.

\begin{lemma}\label{general_quotient_lemma}
Given any $d_0$ and $k_0 < k$ there is a number 
$d_3=d_3(E,d_0,k_0)$ 
such that if $F$ is a general quotient of $E$ of rank $k$ and 
degree $d\ge d_3$, then $F$ has no quotient of degree $d_0$ and rank $k_0$.
\end{lemma}

\noindent
{\bf Definition.} The expression ``general quotient'' means that $F$ is a 
general point of any component of either $\Mg(\Grk,\beta_d)$ or $\Quot$
(and such that the component actually has points corresponding to vector 
bundle quotients).

\medskip
\noindent
{\em Proof of Lemma \ref{general_quotient_lemma}.} If we
count dimensions for the space of the possible $F$'s which {\em do} have 
a quotient $F\rightarrow F_0$ to a bundle $F_0$ of degree $d_0$ and 
rank $k_0$, then 
(viewed as a function of the degree $d$) we have:
\begin{itemize}
\item The dimension of the choice of bundle $F_0$ is constant as a function
of $d$.  Any quotient of such an $F$ is also a quotient of $E$, so the
choice of $F_0$ is at most the dimension of $\Mg(\mathbb{G}(E, k_0)
,\beta_{d_0})$.
\item Fixing an $F_0$ and the kernel $S_{F_0}$ of the surjection 
$E\rightarrow F_0$, we can construct an $F$ by taking a subbundle of $S_{F_0}$
of rank $(r-k)$ and degree $(e-d)$, and then dividing $E$ by this subbundle.
Applying theorem \ref{Mg_dim} to the bundle $S_{F_0}^{*}$ we see that the 
dimension of the choice of subbundles of $S_{F_0}$ is bounded above by a linear function of
$d$ with slope $(r-k_0)$.
\end{itemize}

\noindent
This shows that the dimension of possible $F$'s with the kind of quotient we are 
trying to avoid is bounded above by a linear function
of $d$ with slope $(r-k_0)$.  On the other hand, the dimension of the 
components 
parametrizing vector bundle quotients is bounded below by a linear function of
$d$ with slope $r$ (the cohomological lower bound). Therefore, for large 
enough $d$, no general point of any of these components 
can have a quotient of degree $d_0$ and rank $k_0$ \doneproof.

\medskip
\noindent
{\bf Remark.} This proof shows the utility of Theorems
\ref{Mg_dim} and \ref{Quot_dim}.  
It is important that the theorems put no special condition on
the bundle, since no matter what conditions we put on $E$ we have little control
over $S_{F_0}$, and it is important that the ``slope'' of the upper bound be
the best possible, since it is exactly a small difference in slopes which leads 
to the proof of the lemma.

\medskip

\begin{theorem}\label{there_can_be_only_one}
For all large $d$ there is a unique component of $\Mg(\Grk,\beta_d)$ and 
${\rm Quot}_{k,d}(E)$ whose generic point corresponds to a vector bundle quotient;
this component is of dimension $rd-ke-k(r-k)(g-1)$.
\end{theorem}

\noindent
{\em Proof.} 
The key observation is that there are only a finite number of possible
pairs $(d_0,k_0)$ such that $E$ has a quotient bundle $F_0$ of degree
$d_0$ and rank $k_0$ which also injects into $E\otimes \omega_C$.
If we fix $k_0$, then the degrees of possible quotients of $E$ of rank $k_0$
are bounded from below, and similarly the degrees of possible subbundles
of $E\otimes \omega_C$ are bounded above.  This shows that the number of
possible $d_0$'s for a fixed $k_0$ is finite, and therefore there
are only a finite number of possible pairs, since $k_0$ can only
run between $1$ and $(r-1)$. 

\medskip
\noindent
The reason that this is useful is that if $F$ is any quotient bundle of $E$,
then any map from $F$ into $E\otimes \omega_C$ must factor through one of
the bundles $F_0$ of the type above. If $F$ has no such $F_0$ quotient, 
then the only map from $F$ to $E\otimes \omega_C$ is the zero map, and
this amounts to a vanishing theorem as in Lemma \ref{smooth_point}.

\medskip
\noindent 
Apply Lemma \ref{general_quotient_lemma} to all the possible pairs
$(d_0,k_0)$ above, and pick a $d_3$ which works for all of them.
For any $d\ge d_3$, pick any component of $\Quot$ or $\Mg(\Grk,\beta_d)$
whose general point corresponds to a vector bundle quotient.
By Lemma \ref{general_quotient_lemma}, the general point $F$ of this 
component has no map to $E\otimes \omega_C$, and so (by Serre duality) both
$h^1(E^{*}\otimes F)$ and $h^1(S_F^{*}\otimes F)$ are zero.

\medskip 
\noindent
Pick a deformation (as an abstract bundle) of $F$ over an irreducible base 
$B$ so that the generic
member of the deformation is a stable bundle.  
Throw away the closed set of points of $B$ 
where either $h^1(E^{*}\otimes F_b)$ or $h^1(S_{F_b}^{*}\otimes F_b)$ jumps.  
We now have $h^1(E^{*}\otimes F_b)=0$
for all $b$ in $B$, so that the vector spaces $\Hom(E,F_b)$ form a bundle
over $B$.  The condition of being a surjection is open, so that we can also
deform the surjection $E\rightarrow F$ to a surjection $E\rightarrow F_b$ for
any general point $b$ of $B$. The deformation was constructed so that for 
general $b$, the bundle $F_b$ is stable.

\medskip
\noindent
Since $h^1(S_{F_b}\otimes F_b)=0$ for all $b$ in $B$, this deformation 
always stays in the smooth locus of $\Mg(\Grk,\beta_d)$ or $\Quot$, therefore
in deforming this surjection we remain on the same component that we started on.

\medskip
\noindent
This shows that as long as $d\ge d_3$, the condition that the general point
of a component correspond to a vector bundle quotient is enough to imply that
the general quotient on that component is actually stable.  From
Proposition \ref{good_component} we know that
once $d$ is large enough, there is a unique such component, of dimension
$rd-ke-k(r-k)(g-1)$, and this proves the theorem \doneproof.

\medskip
\noindent
{\bf Remark.} As the referee points out, the last part of the proof above can be 
simply summarized as follows: we have that the Kodaira-Spencer
infinitesimal deformation morphism at a point of the corresponding 
component chosen as above is surjective, and thus a generic 
quotient has to be stable.

\medskip
\noindent
{\bf Remark.} The condition $d\ge d_3$ is already enough to guarantee 
the uniqueness of the component parametrizing vector bundle quotients, since
(as the proof shows) every such component has quotients which form an open
set of the set of rank $k$, degree $d$, stable bundles.  Since this set is
irreducible, and these quotients lie in the smooth loci of the components,
the uniqueness statement follows.

\medskip
\noindent
The theorem immediately implies the following useful Corollary:

\medskip

\begin{corollary}
Every vector bundle $E$ on $C$ can be written as an extension
$$0\longrightarrow E_1\longrightarrow E\longrightarrow 
E_2\longrightarrow 0,$$
with $E_1$ and $E_2$ generic stable bundles of prescribed ranks.
\end{corollary}

\medskip
\noindent
Note that we are free to choose what we mean by generic, that is, we may
specify ahead of time any open condition on the space of stable bundles 
(for example 5.4 stable) and then express $E$ as an extension where 
both $E_1$ and $E_2$ satisfy that open condition.  Similarly, we can freely specify
the ranks $r_1$, $r_2$ of the bundles $E_1$, $E_2$ involved in the extension.

\bigskip
\noindent
We now study when the spaces $\Quot$ or $\Mg(\Grk,\beta_d)$ can be irreducible
for large $d$.

\medskip
\noindent
Our main result says that in the case of the Quot scheme the best possible
result is true:

\begin{theorem}\label{irred_of_Quot}
For any vector bundle $E$ on $C$, there is an integer $d_Q = d_Q(E,k)$ such that
for all $d\geq d_Q$, ${\rm Quot}_{k,d}(E)$ is irreducible. 
\end{theorem}

\medskip
\noindent
{\bf Remark.} Of course, this unique component is generically smooth,
of dimension $rd- ke - k(r-k)(g-1)$, and the 
general point corresponds to a vector bundle quotient,
all of which follows immediately from Theorem \ref{there_can_be_only_one}.

\medskip
\noindent
{\em Proof of Theorem \ref{irred_of_Quot}.}
For simplicity, let us use ${\rm Quot}_{k,d}^0 (E)$ for the open subset 
of $\Quot$ parametrizing vector bundle quotients 
(note that for some $d$ this might miss entire components of $\Quot$). 

\medskip
\noindent
Also set: 
$$w_d = \dim {\rm Quot}_{k,d}^0(E)  - \left(dr-ke-k(r-k)(g-1)\right),$$
i.e. the difference between the dimension of this subset and the
expected dimension (the lower bound).
We have $w_d \ge 0$ for all $d\ge d_k$, and the $w_d$ are eventually zero
(by Theorem 6.2) for large $d$.

\medskip
\noindent
Let $W$ be the set of $d\geq d_k$ so that $w_d\ne0$, and define $M$ by:
$$M := \max_{d\in W} \left\{d+\frac{w_d}{k}\right\}.$$
Finally, set 
$$d_Q = \max(\lceil M\rceil, d_3(E,k)),$$
where $\lceil M\rceil$ is the smallest integer greater than $M$ and $d_3(E,k)$ is the 
integer $d_3$ appearing in the proof of Theorem \ref{there_can_be_only_one} 
(i.e., such 
that there is a unique component corresponding to vector bundle quotients for $d\geq d_3$).

\medskip
\noindent
The claim is that for all $d> M$, the general point of every component
of $\Quot$ corresponds to a vector bundle quotient.  Since, by Theorem 6.2 there
is only one of these as soon as $d\ge d_3$, this proves the irreduciblity
of the Quot scheme for all $d\ge d_Q$.

\medskip
\noindent
Suppose we pick a degree $d'$ smaller than $d$ and construct, as in Proposition \ref{21},
torsion quotients of degree $d$ by introducing torsion of order $\delta=(d-d')$ 
(all torsion loci are produced by this kind of construction).
The torsion quotients we get this way have dimension
$$\dim \mbox{Quot}_{k,d'}^0(E) + \delta(r-k)$$
or, another way to write this (using the $w$'s) is
$$w_{d'} + d'r -ke -k(r-k)(g-1) + \delta(r-k).$$
The difference between this dimension and the lower 
bound $rd -ke - k(r-k)(g-1)$ is:
$$-(d-d^{\prime})r + w_{d^{\prime}} +\delta (r-k) =  w_{d^{\prime}} -
\delta k.$$
But, by assumption, $d > M$, which means that $d > d^{\prime} 
+\frac{w_{d^{\prime}}}{k}$,
or equivalently $w_{d^{\prime}} - \delta k < 0$ (or, if $d^{\prime}\notin
W$, then $w_{d^{\prime}} =0$
and what we want is automatic). In other words, the dimension of the space of quotients 
with torsion is strictly lower than
the lower bound, so the general point must be an actual vector bundle quotient,
which finishes the proof \doneproof.

\medskip
\noindent
{\bf Remark.} We will use Theorem \ref{irred_of_Quot} at the end of 
the next section to show that $\Mg(\Grk,\beta_d)$ is connected 
for large $d$ (see Proposition \ref{connected_Mg}).

\medskip
\noindent
{\bf Remark on $g=0$.} In the genus zero case, when $E$ is the trivial bundle, Str{\o}mme 
\cite{Stromme} has shown that the Quot scheme is smooth and irreducible for
every $k$ and $d$.  For a general vector bundle $E$ over $\PP^1$, neither statement
need be true.  However, the asymptotic results of 
Proposition \ref{good_component} and Theorems
\ref{there_can_be_only_one} and \ref{irred_of_Quot} still remain valid for genus zero,
with a slight modification in the proof.  The modification is necessary because
there are not many semistable bundles on genus $0$ curves (if $d$ is not a multiple of
$r$, $E$ cannot be semistable) and we use the semistable bundles as a means of
connecting and analyzing the components.  To replace them, for each $d$ and each $k\geq 1$, 
let $F_{k,d}$ be the unique (up to isomorphism) bundle of rank $k$ and degree $d$ on $C\cong\PP^1$ 
such that $F_{k,d}= \oplus_i\Osh_C(a_i)$ with $|a_i-a_j|\leq1$ for all $1\leq i,j\leq k$.  Then the results
in the proposition and theorems listed above are true in genus zero, the only changes being that the 
words ``stable bundle'' used in the statements and arguments should be changed to
``of type $F_{k,d}$''.

\medskip
\noindent
We now study the irreducibility of the spaces $\Mg(\Grk,\beta_d)$.
Unlike the case of $\Quot$,
it is typical (even for large $d$) that there are components of larger 
dimension consisting entirely of reducible maps.
If we want to ensure that $\Mg(\Grk,\beta_d)$ has no other components,
then we must put strong conditions on the vector bundle $E$.
(We assume in what follows that $g\geq 2$, see the remark after the proof for
a discussion of $g=0,1$.)

\begin{proposition}
\label{irred_comps}
If $\Mg(\Grk,\beta_d)$ 
is irreducible for all $k$ and all large $d$, then
$E$ is a stable bundle, generic in the sense of Example 
5.\ref{generic_stable_bundle}. 
\end{proposition}

\medskip
\noindent
{\em Proof of Proposition \ref{irred_comps}.}
Let $E$ be a vector bundle of rank $r$ and degree $e$ satisfying
the conditions of the proposition. Consider $k$ fixed; we will show that 
(at least as far as rank $k$ quotients are concerned) $E$ is a generic 
stable bundle in the sense of Example 5.4. Since it holds for all $k$
this will prove the proposition.

\medskip
\noindent
For $d$ large we know that there is a component of dimension 
$rd-ke-k(r-k)(g-1)$, and since we are assuming $\Mg(\Grk,\beta_d)$ irreducible,
this is the only component. Let $m_k$ be (as in \S3) the dimension of 
$\Mg(\Grk,\beta_{d_k})$. 

\medskip
\noindent
In Proposition \ref{20} we constructed a family of reducible
curves of dimension $m_k + (d-d_k)r -1$.  This family has to lie in the unique
component, and since the curves constructed have generically only one node,
it is a divisor in this component, which gives us 
$$ m_k + (d-d_k)r = rd-ke-k(r-k)(g-1),\,\,\mbox{or}$$
$$m_k = d_kr -ke -k(r-k)(g-1).$$
Since $m_k$ is nonnegative, the right hand side of the equation is as well.
(The fact that the subset above is a divisor is not really necessary
for the argument -- we only need a proper subset, as the reader can
easily check.) 

\medskip
\noindent
For the fixed numerical data of $r$, $e$, $g$, and $k$, let $D_k$ be 
the smallest integer so that $D_kr-ke-k(r-k)(g-1)$ is nonnegative.  By the 
results in \cite{Teixidor} this is
the smallest degree in which a generic stable bundle of rank $r$, 
degree $e$ has a rank $k$ quotient.

\medskip
\noindent
The fact that $m_k$ is nonnegative gives us $d_k \ge D_k$, but since the 
number $d_k(E)$ is {\em lower} semicontinuous in families of vector bundles, 
and since any bundle can be deformed to a generic stable bundle, 
we always have $d_k(E)\le D_k$, and this gives $d_k=D_k$.

\medskip
\noindent
This also shows that $E$ has no rank $k$ destabilizing quotient bundles, either
from the fact that $(d_kr-ke)$ is positive, or from the
fact that the degrees of its rank $k$ quotients are the same as that of a
generic stable bundle of the same rank and degree.

\medskip
\noindent
To see that $E$ is 5.4 stable (that is, has the same properties listed in 
example 5.4) we need to show the following additional two things:
\begin{enumerate}
\item $\Mg(\Grk,\beta_d)$ has a component of dimension $dr-ke-k(r-k)(g-1)$ 
for all $d\ge d_k$.
\item There are no components of larger dimension.
\end{enumerate}

\medskip
\noindent
{\bf Remark.} We really only need to show this for ``small $d$'', since for
large $d$ this is a consequence of the assumptions and of Proposition \ref{good_component}.

\medskip
\noindent
The families of nodal curves constructed in Proposition \ref{20} are of dimension
$m_k + (d-d_k)r -1$ which is $dr-ke-k(r-k)(g-1) -1$ in our case. 

\medskip
\noindent
The deformation theoretic lower bound for this family of curves is 
$dr-ke-k(r-k)(g-1)$, which is one larger.  Therefore 
$\Mg(\Grk,\beta_d)$ has a component of dimension at least
$dr-ke-k(r-k)(g-1)$ for all $d\ge d_k$. We will be finished once we 
establish condition $2$.

\medskip
\noindent
Suppose for some $d_1$ that there is a component of $\Mg(\Grk,\beta_{d_1})$ 
of dimension strictly greater than 
$d_1r-ke-k(r-k)(g-1)$.  Again by attaching a single rational curve 
mapping in with large degree 
we can produce nodal families of dimension at least $dr-ke-k(r-k)(g-1)$ for all $d\ge d_1$.  

\medskip
\noindent
For large $d$, this is 
greater than or equal to the dimension of the unique component.  
By \ref{there_can_be_only_one} there is also a component of
this dimension where the general member of the universal family is a smooth
curve.  This contradicts the uniqueness assumption and gives condtion 2 \doneproof.

\medskip
\noindent
The fact that the components are always of the same dimension
as the deformation theoretic lower bound has an important consequence. 
Combined with Theorem \ref{there_can_be_only_one} this consequence proves the 
converse of Proposition \ref{irred_comps}.
 
\begin{proposition}\label{always_quotient}  
For a stable vector bundle $E$, generic in the sense
of example 5.\ref{generic_stable_bundle}, 
the generic point of every component of $\Mg(\Grk,\beta_d)$ corresponds to 
a map from a smooth curve $($i.e., vector bundle quotient$)$. 
\end{proposition}

\noindent
{\em Proof.} Suppose that we have a component consisting
entirely of reducible maps.  We consider this to be built up (as usual) from
a smaller component by attaching rational curves.  We know the dimension of
the smaller component, and if we compute the dimension of the family we obtain
by gluing on the curves, we find that it is of {\em smaller} dimension than
it is supposed to be \doneproof.

\bigskip
\noindent
Putting together Propositions \ref{irred_comps} and \ref{always_quotient}, and 
Theorem \ref{there_can_be_only_one}, we arrive at

\begin{theorem}\label{irred_comps_thm}
The spaces $\Mg(\Grk,\beta_d)$ are irreducible for all $k$ and
all large $d$ if and only if $E$ is stable bundle, generic in the sense of 
Example 5.\ref{generic_stable_bundle}. 
\end{theorem}

\noindent
{\bf Remark.}  In the proof of Proposition \ref{irred_comps} we showed that
if for a fixed $k$, $\Mg(\Grk,\beta_d)$ is irreducible for all large $d$, then
$E$ is 5.4 stable ``in rank $k$'' (i.e. at least as far as rank $k$ quotients 
are concerned).  On the other hand, if we assume that $E$ is 5.4 stable in 
rank $k$,
then the result of Proposition \ref{always_quotient} is still true for rank $k$
quotients, and so applying Theorem \ref{there_can_be_only_one} we get a slight
generalization of Theorem \ref{irred_comps_thm}:  For fixed $k$, 
$\Mg(\Grk,\beta_d)$ is irreducible for all large $d$ if and only
if $E$ is 5.4 stable in rank $k$.

\medskip
\noindent
{\bf Remark.} In the case that $g=0$ or $g=1$, the proof breaks down at the
following point: when we wanted to show that $D_k \ge d_k(E)$ we used the
fact that $d_k$ is lower semicontinuous in families, and that we could deform
$E$ to a more general bundle $E'$ with $d_k(E')=D_k$.  In the case that $g<2$,
the gap is that we don't know that such a bundle $E'$ exists.  The rest of
the argument still goes through, and shows that $\Mg(\Grk,\beta_d)$ is 
irreducible for large $d$ if and only if all components of $\Mg(\Grk,\beta_d)$
have exactly the dimension of the lower bound, $dr-ke-k(r-k)(g-1)$.  The only
difference in the result is that we don't know when $E$ starts to have quotients
of rank $k$, i.e., we don't know what the value of $d_k(E)$ is -- it could
be strictly larger than $D_k$.  As an illustration, in Example 
5.\ref{trivial_example} 
in the case that $g=0$ the value of $D_k$ is negative,
but $d_k(E) =0$ for all $k$.  However, the dimension of $\Mg(\Grk,\beta_d)$
is the predicted value for all $d\ge 0$, and these spaces are easily seen
to be irreducible.

\medskip
\noindent
Since in the case of 5.4 stable bundles one has good numerical control over the invariants
involved, all the steps in the proof can be made effective, and we have:
                                                                                                                             
\begin{theorem}\label{effective_bound}
Let $E$ be a 5.4 stable bundle. Then, for a fixed $k$, $\Mg(\Grk,\beta_d)$ 
and ${\rm Quot}_{k,d}(E)$ are irreducible for
$$d\geq \frac{ke + k(r-k)(g-1)}{r} + 2(k-1)(g-1) + (r-k)(k-1).$$
In particular
$$d\geq (r-1)\mu(E) + \left(\frac{9}{4} r - 4\right)(g-1) + \frac{r^2}{4} - r +1$$
works for all $k$.
\end{theorem}
 
\noindent
{\em Proof.}
As noted after the proof of Theorem \ref{there_can_be_only_one}, all we need to do is to find
an effective expression for the number $d_3$ appearing in Lemma \ref{general_quotient_lemma}
which works for all the possible pairs $(k_0, d_0)$ which are the rank and degree of quotients
of $E$ that are also subsheaves of $E\otimes \omega_C$ and such that $k_0< k$.
 
\medskip
\noindent
Fixing such a pair $(k_0, d_0)$,
we first bound the dimension of the space of quotients $F$ of $E$ which in turn have
a quotient $F_0$ of rank $k_0$ and degree $d_0$. These are themselves quotients of $E$,
and since $E$ is 5.4 stable, we have
$${\rm dim}\{F_0\} \leq {\rm dim~}{\rm Quot}_{k_0,d_0}(E) = rd_0 - k_0 e - k_0(r-k_0)(g-1).$$
Now each such $F_0$ determines an exact sequence
$$0\longrightarrow S_{F_0}\longrightarrow E\longrightarrow F_0\longrightarrow 0,$$
and so the space of $F$'s which have $F_0$ as a quotient is determined by the space of
subbundles of $S_{F_0}$ of rank $r-k$ and degree $e-d$. Its dimension is bounded above by
the formula given in Theorem \ref{Quot_dim}, and after a small computation we get that the dimension
of the space of all $F$'s we are looking at is bounded by:
$${\rm dim}\{F\}\leq rd_0 - k_0 e - k_0(r-k_0)(g-1) + (r-k)(k-k_0) + (d - d_k(E))(r-k_0).$$
On the other hand (since $E$ is 5.4 stable) the dimension of degree $d$ rank $k$
quotients is exactly $rd - ke - k(r-k)(g-1)$, and so
we want $d$ large enough so that
\begin{eqnarray*}
k_0d  &>&  ke + k(r-k)(g-1) + rd_0 - k_0 e \\
 & & \mbox{} - k_0(r-k_0)(g-1) - d_k(E)(r-k_0) + (r-k)(k-k_0).
\end{eqnarray*}
 
\medskip
\noindent
Again using the fact that $E$ is 5.4 stable, we have
$$-d_k(E)\leq \frac{-ke - k(r-k)(g-1)}{r}.$$
Plugging this into the previous inequality and computing further, we obtain
$$d >  \frac{ke + k(r-k)(g-1)}{r} +r\mu(F_0) - e - (r-k_0)(g-1) +
(r-k)\left(\frac{k}{k_0} - 1\right).$$
In the proof of Theorem \ref{there_can_be_only_one} we are assuming that $F_0$
is a subsheaf of $E\otimes \omega_C$ (which is also 5.4 stable), giving the
additional condition:
$$\mu(F_0) \leq \mu(E)  + \frac{r+k_0}{r}(g-1).$$
 
\medskip
\noindent
And this gives us:
$$ d > \frac{ke + k(r-k)(g-1)}{r} + 2k_0(g-1) +
(r-k)\left(\frac{k}{k_0} - 1\right).$$
Using the inequalities $1\leq k_0\leq k-1$, we see that this is satisfied
whenever:
$$d \ge \frac{ke + k(r-k)(g-1)}{r} + 2(k-1)(g-1) + (r-k)(k-1),$$
which is the desired bound on $d$ depending only on $k$, $r$ and $e$.
 
\medskip
\noindent
Optimizing the terms involving $k$ individually, we arrive at the final answer which works
for all $k$: all spaces $\Mg(\Grk,\beta_d)$ and ${\rm Quot}_{k,d}(E)$ are irreducible for:
$$d\geq d_3 = (r-1)\mu(E) + \left(\frac{9}{4} r - 4\right)(g-1) + \frac{r^2}{4} - r +1 \,\,\doneproof.$$
 
\medskip
\noindent  
{\bf Remark.} In the final two steps of the computation above we optimized each term containing
$k_0$ and then $k$ individually,
rather than optimizing the whole expression as a function of $k$.
Although the bound obtained is 
weaker, we chose to do this in order to avoid some extra, not very illuminating, calculations.
The reader can certainly find better approximations for any particular example, if needed.

\medskip
\noindent
In the case when $E$ is only assumed to be stable, the number $d_3$ can still be approximated 
(thus making Theorem \ref{there_can_be_only_one} effective), by simply changing the 5.4 
stability condition into usual stability. 
We skip the computation, which is identical to the one above, but list the result:
$$d \geq (r-1)\mu(E) + \left(\frac{r^2 + 8r}{4}\right)(g-1) + \frac{r^2}{4} + r.$$

\bigskip
\noindent
It is worth noting that, at least in characteristic $0$, 
we only need $E$ to be stable in order to perform some other effective 
calculations. We conclude with a sample such result, which essentially makes Proposition
\ref{good_component} effective: 

\begin{proposition}
Assume that $C$ is defined over a field of characteristic $0$, and let $E$ be a stable bundle of 
rank $r$ and degree $e$ on $C$. Then, for a fixed $k$, all semistable bundles of rank $k$ and
degree $d\geq k(\mu (E) + 2g)$ are quotients of $E$, and in fact smooth points of a unique 
component of ${\rm Quot}_{k,d}(E)$. In particular
$$d\geq (r-1)(\mu(E) + 2g)$$
works for all $k$.  
\end{proposition}

\noindent
{\em Proof.}
Following the proof of Proposition \ref{good_component} and the subsequent Remark, we see that we only need
to bound the numbers $d_1$ and $d_2$ appearing in Lemmas \ref{good_homs} and \ref{smooth_point} respectively.
This involves only the well known facts that the tensor product of two semistable bundles is semistable
(in characteristic $0$) and that for a semistable bundle $V$, $\mu(V) \geq 2g-1$ implies $h^1(V)=0$ and 
$\mu(V)\geq 2g$ implies global generation.

\medskip
\noindent
In Lemma \ref{good_component} we had $\mu(E^* \otimes F) = \frac{d}{k} - \mu(E)$, and since $E$ 
and $F$ are now both stable, $E^*\otimes F$ will be globally generated and have vanishing $H^1$ as 
long as $d\geq k(\mu(E) + 2g)$.
Even simpler, the definition of stability for $E$ implies that the sequence of  
inclusions in Lemma \ref{smooth_point} cannot hold if $d\geq k(\mu(E) + 2g-2)$. This gives the desired bound    
\doneproof.

\newpage
\noindent
\mysection{Maps from $\Mg(\Grk,\beta_d)$ to $\Quot$}
\label{maps_section}

\medskip
\noindent
We would like to try and extend the isomorphism between the open sets of
$\Mg(\Grk,\beta_d)$ and $\Quot$ parametrizing vector bundle quotients to
a morphism between the two spaces.  Since the boundary loci of the stable
map space are always of larger dimension than the boundary loci of
the Quot scheme, the only possible direction for such a morphism is from 
$\Mg(\Grk,\beta_d)$ to $\Quot$.

\medskip
\noindent 
The moduli functor for $\Mg(\Grk,\beta_d)$ classifies flat families 
$\Cfam\rightarrow B$ of nodal curves over a base $B$ along with a map
$f:\Cfam \longrightarrow \Grk \times B$ such that over each point $b$ of $B$,
the fibre $\Cfam_b$ and morphism $f_b$ is a stable map into $\Grk$, with
$f_{b*}[\Cfam_b]=\beta_d$.  

\medskip
\noindent
The corresponding moduli functor for $\Quot$ classifies flat families $\Qfam$ 
of quotients of $E$ on $C \times B$, such that over each point $b$ the 
restriction $\Qfam_b$ of the family to the fibre $C$ over $b$ is a coherent
sheaf of degree $d$ and generic rank $k$. (We will also use the symbol
$E$ for the pullback $p_1^{*}E$ on $C\times B$).

\medskip
\noindent
To give a morphism from $\Mg(\Grk, \beta_d)$ to $\Quot$ is to give a 
functorial way of producing a flat family of quotients $\Qfam$ on 
$C \times B$ from a family of stable maps over $B$.  

\medskip
\noindent
The Grassmann bundle $\Grk$ comes with the projection 
$\pi:\Grk \longrightarrow C$, and so given any family of stable maps
over $B$, we can compose the morphism $f:\Cfam \longrightarrow \Grk \times B$ 
with $\pi$ to get a map 
$$g=(\pi\circ f):\Cfam \longrightarrow C\times B.$$

\medskip
\noindent
One way to produce sheaves on $C\times B$ is to push forward bundles from the
family $\Cfam$, and we will use the following lemma (with $g=(\pi\circ f)$ 
as above).

\begin{lemma}\label{base_change_lemma}
If $\Vbndle$ is a vector bundle on a family of stable maps over $B$ 
such that $R^1g_{*}\Vbndle=0$, then $g_{*}\Vbndle$ is a flat family of coherent
sheaves on $C\times B$, and formation of $g_{*}\Vbndle$ commutes with
arbitrary base change $B^{\prime}\rightarrow B$. 
\end{lemma}

\noindent
{\em Proof.} This is a relative version of the standard theorem on base change
and cohomology, and the proof is almost the same. 
The bundle $\Vbndle$ is flat over $B$ since $\Cfam$ is a flat family 
and $\Vbndle$ locally free.   We can compute the derived sheaves on $C\times B$
by means of the relative \v{C}ech complex, which is a two term complex 
$\Vbndle^{0}\rightarrow \Vbndle^{1}$ of quasi-coherent flat $\Osh_B$ modules.  The
map of complexes is surjective since by hypothesis $R^1g_{*}\Vbndle=0$. 
This gives us the exact sequence 
$$ 0\longrightarrow g_{*}\Vbndle \longrightarrow \Vbndle^{0} \longrightarrow
  \Vbndle^{1} \longrightarrow 0 $$
on $C\times B$.  Since both $\Vbndle^{0}$ and $\Vbndle^{1}$ are flat over $B$, 
$g_{*}\Vbndle$ is as well. Computation of the kernel commutes with
arbitrary base change induced from $B$, again since $\Vbndle^{1}$ is flat, and
finally, $g_{*}\Vbndle$ is coherent over $C\times B$ since $g$ 
is proper \doneproof.

\medskip
\noindent
The map $g$ from $\Cfam$ to $C\times B$ is the contraction of all the 
``rational tails'' on the curves, and the derived sheaf $R^1g_{*}\Vbndle$ 
vanishes whenever the bundle $\Vbndle$ has no $H^1$ when restricted to these 
fibres.  If it happens that $\Vbndle$ is trivial when restricted to the fibres,
then the pushforward $g_{*}\Vbndle$ is a vector bundle on $C\times B$, a fact
which we will use later on.

\medskip
\noindent
Given any family of stable maps, the exact sequence of tautological bundles 
$$0\longrightarrow \Sbndle \longrightarrow \pi^{*}E\longrightarrow 
\Fbndle \longrightarrow 0$$
on $\Grk$ pulls back via $f$ to give a similar exact 
sequence on the family $\Cfam$.  For any rational curve $D$ mapping to a fibre
of the relative Grassmannian, the pullback $f^{*}\pi^{*}E$ is trivial; this
means that the tautological subbundle $f^{*}\Sbndle$ restricted to $D$ 
decomposes as $\oplus_i \Osh(a_i)$ with all $a_i\le 0$, and 
similarly that the pullback of
the tautological quotient bundle $f^{*}\Fbndle$ restricted to $D$ decomposes
as a sum $\oplus_j \Osh(b_j)$ with all $b_j\ge 0$.

\medskip
\noindent
The decomposition of $f^{*}\Fbndle$ on the rational curves shows that this
bundle has no $H^1$ on any of the fibres of the map $g:\Cfam
\rightarrow C\times B$,
and therefore Lemma \ref{base_change_lemma} says that the sheaf
$g_{*}(f^{*}\Fbndle)$ forms a flat (over $B$) family of coherent sheaves on
$C\times B$.  On each fibre $C_b$, the restriction is of degree $d$ and 
rank $k$.

\medskip
\noindent
This family of coherent sheaves unfortunately does {\em not} in 
general give a map from 
$\Mg(\Grk,\beta_d)$ to $\Quot$.   There are two ways to see this;
the first is that when the subbundle $f^{*}\Sbndle$ has higher cohomology
on the rational fibres, the induced map from $g_{*}(f^{*}\pi^{*}E)=E$
to $g_{*}(f^{*}\Fbndle)$ on $C\times B$ is not surjective.  
This does not completely rule out this family, since it is conceivable
that it might be a quotient of $E$ in some other way.  What really shows
that this is the wrong construction is the following.
If $C_b$ is a fibre where the family $g_{*}(f^{*}\Fbndle)$ has torsion, and
$x$ a local coordinate at a torsion point, then multiplication by $x$ kills
all torsion at that point.  In general for quotients of $E$ with torsion,
we expect situations where multiplication by $x^2$ (or higher powers of $x$)
is necessary to annihilate the torsion submodule.  

\medskip
\noindent
In situations where the subbundle does {\em not} have higher cohomology along 
the rational fibres
(for instance, if $d=d_k+1$, where the $a_i$'s are all zero except for one
$-1$), this does give the right morphism to $\Quot$. 

\medskip
\noindent
The correct approach is to try and alter the subbundle $f^{*}\Sbndle$ along
the contracted rational fibres of $g$ so that it 
is trivial along these fibres.  The pushforward
will then be a vector bundle on $C\times B$ which includes into $E$, 
and gives a flat family of quotients.

\medskip
\noindent
This is always possible in the case that $k=(r-1)$, when the subbundle is a line
bundle.

\begin{theorem}\label{(r-1)_map}
If $k=(r-1)$ then there is a surjective morphism from $\Mg(\Grk,\beta_d)$ to
${\rm Quot}_{k,d}(E)$ which extends the map on the locus where the domain curve is smooth.
\end{theorem}

\medskip
\noindent
{\em Proof.}  In order to describe a construction to be applied naturally 
(i.e., functorially) to all families of stable maps, we should describe
the
construction on $\Mgi(\Grk,\beta_d)$, the universal curve over 
$\Mg(\Grk,\beta_d)$.  The idea is that the loci of rational tails form divisors
in $\Mgi(\Grk,\beta_d)$ and twisting by these divisors will allow us to 
make the line bundle $f^{*}\Sbndle$ trivial along all the tails.  
There is some awkwardness in the description because there can be
components of $\Mg(\Grk,\beta_d)$ where the general map is reducible, and 
here the loci of rational tails we must twist by in $\Mgi(\Grk,\beta_d)$
occur as both divisors and components. 

\medskip
\noindent
Let $D$ be either a divisor or a component of $\Mgi(\Grk,\beta_d)$ such
that the general fibre of the map to $\Mg(\Grk,\beta_d)$ is an irreducible
rational curve, with the additional condition that this general fibre touches
only one other curve in the domain of the map it is 
associated to.  (This extra condition only comes up because of the components 
where the general map is reducible.)

\medskip
\noindent
For each such locus $D$ there is a well defined line bundle $\Osh_{\Mgi}(D)$ on
$\Mgi(\Grk,\beta_d)$\footnote{
The fact that O(D) is a line bundle can be seen by pulling
it back from the (nonseparated) stack of all nodal curves of 
genus g, where deformation theory shows it to be a line bundle.  This argument 
is well known to the experts, but there does not seem to be any
published reference.};
this line bundle has degree $-1$ when restricted to the general fibre of the
induced map from $D$ to $\Mg(\Grk,\beta_d)$.  
We attach the weight
$\delta$ to such a $D$ if the degree of $f^{*}\Sbndle$ restricted to the
general fibre in $D$ is $-\delta$.  

\medskip
\noindent
The fact that 

\medskip
\noindent
Set
$$ \Sbndle^{\prime}=f^{*}\Sbndle
\otimes_{\delta}
\Osh_{\Mgi}(-\delta \cdot D_{\delta}),$$
where the tensor product runs over all possible $\delta$ and all possible
$D_{\delta}$ with weight $\delta$ satisfying the conditions above.

\medskip
\noindent
For any family $\Cfam$ of stable maps over a base $B$, this construction
gives us a line bundle $\Sbndle^{\prime}$ on $\Cfam$ of degree $0$ on 
all components of all rational tails.  In addition, this bundle has
a nonzero map $\Sbndle^{\prime}\rightarrow f^{*}\Sbndle$ obtained by 
multiplying by a section of $\otimes_{\delta} \Osh(\delta\cdot D_{\delta})$. 

\medskip
\noindent
By Lemma \ref{base_change_lemma} the pushforward $g_{*}\Sbndle^{\prime}$
is a flat family of line bundles on $C\times B$, and comes with a natural 
inclusion $g_{*}\Sbndle^{\prime} \hookrightarrow E$ induced from the map
$\Sbndle^{\prime}\rightarrow f^{*}\Sbndle$ on the family $\Cfam$.  Let
$\Qfam$ be the quotient.

\medskip
\noindent
By definition, $\Qfam$ is a family of quotients of $E$.  
If $b$ is a point of the base $B$ where the domain of the map $\Cfam_b$
is an irreducible curve, then the construction with $\Sbndle^{\prime}$
changes nothing near $\Cfam_b$ and the quotient on the fibre $C_b$ is
precisely the quotient we expect.  The only step remaining is to see that
$\Qfam$ is a flat family over $B$.

\medskip
\noindent
Let $\Id$ be an ideal sheaf in $\Osh_B$.  Tensoring the defining sequence
for $\mathcal{Q}$ with $(\Osh_B/\Id)$ gives the sequence
$$0\rightarrow \mathcal{T}or_1(\Osh_B/\Id,\Qfam)\longrightarrow g_{*}\Sbndle^{\prime}
\otimes (\Osh_B/\Id) \longrightarrow
E \otimes (\Osh_B/\Id)
\longrightarrow \mathcal{Q} \otimes (\Osh_B/\Id) \rightarrow 0.$$
However, both $g_{*}\Sbndle^{\prime}$ and $E$ are vector bundles, and at
a general point of every fibre over $B$ the map $g_{*}\Sbndle^{\prime}
\rightarrow E$ is an injection, therefore the map from $g_{*}\Sbndle^{\prime}
\otimes (\Osh_B/\Id)$ to $E\otimes (\Osh_B/\Id)$ is
an injection and $\mathcal{T}or_1(\Osh_B/\Id,\Qfam)=0$ \doneproof.

\medskip
\noindent
It is easy to see what this construction is doing.  Suppose that we start
with $C_0$, 
a section of the relative Grassmannian $\Gr(E,r-1)$, and let $S_0$ be the
universal sub-line bundle restricted to $C_0$.  Attach 
a rational curve to $C_0$ at a point $p$,  and extend the
inclusion $C_0 \hookrightarrow \Gr(E,r-1)$ to this reducible curve by mapping 
the rational component in so that it accounts for $\delta$ in degree. 
In this situation the
construction produces the subsheaf $S^{\prime}$ of $S_0$
(on $C_0$) which vanishes to order $\delta$ at $p$.  This bundle includes 
into $E$, and gives a quotient with torsion of order $\delta$ at that point.

\medskip
\noindent
This suggests how the ``right'' torsion quotient to associate to any 
stable map into the relative Grassmannian might look, in the general case
(i.e., for arbitrary $k$).

\medskip
\noindent
{\em Possible description of Quotient:} Suppose that we have a reducible
map consisting of a section $C_0$ and a single rational curve attached
to $C_0$ at $p$ mapping into $\Grk$.  Let $S_0$ be the restriction of the
universal subbundle to the section $C_0$ and let $S_1$ be the restriction
to the rational curve.  The bundle $S_1$ on the rational curve breaks up
into a sum $\oplus_{i=1}^{r-k} \Osh(a_i)$ with $a_i \le 0$.

\medskip
\noindent
We can identify the fibre of $S_0$ at $p$ with the fibre of $S_1$ at the 
corresponding point of attachment.  The picture in the case $k=(r-1)$ suggests
that we should be able to choose a local basis 
$\{s_1,\ldots,s_{(r-k)}\}$ for $S_0$ at $p$, compatible with the splitting
on $S_1$, so that 
if $x$ is a local coordinate at $p$ on $C_0$, the subsheaf $S^{\prime}$ of 
$S_0$ generated locally by $\{x^{-a_1}s_1,\ldots, x^{-a_{(r-k)}}s_{(r-k)}\}$
gives the correct quotient.

\medskip
\noindent
This unfortunately is also not true in most situations (see Example
\arabic{section}.\ref{strange_limit} below).  It is however
true in two special cases, and these special cases will be sufficient to show
that there is in general {\em no} map from $\Mg(\Grk,\beta_d)$ to $\Quot$
(if $k\ne (r-1)$) which extends the map on smooth curves.

\begin{proposition}\label{generic_splitting_case}
Suppose that $\Cfam\longrightarrow B$ is a one parameter family of stable maps 
with both $\Cfam$ and $B$ smooth,
and that at a point $b_0$ of $B$ the fibre
$\Cfam_{b_0}$ consists of a section $C_0$ and a single rational curve joined
at a point $p$. If the pullback $f^{*}\Sbndle$ of the universal subbundle
decomposes on the rational curve as $\oplus \Osh(a_i)$ such that 
$|a_i-a_j|\le 1$ for all $i$, $j$, then the limiting quotient $\Qfam_{b_0}$ 
in ${\rm Quot}_{k,d}(E)$ is given by the rule above.
\end{proposition}

\noindent
{\bf Remark.} The condition that both $\Cfam$ and $B$ are smooth means that
the general fibre is as well, and this gives a generically defined map from
$B$ to $\Quot$; it is this map that we use to compute the one parameter limit.

\medskip
\noindent
{\bf Remark.} In this case (that $|a_i-a_j|\le 1$) {\em any} choice of local
generators of the sheaf $S_0$ at $p$ (compatible with the splitting along
the rational curve) gives the same subbundle $S^{\prime}$. (See Proposition
\ref{one_(-2)_case} for an example where this choice does matter).

\medskip
\noindent
{\em Proof.} Let $D$ be the rational curve over $b_0$ on the smooth surface 
$\Cfam$.  Suppose that the subbundle $f^{*}\Sbndle$ restricted to $D$ 
decomposes as 
$$f^{*}\Sbndle|_{D} =\, 
\stackrel{m}{\oplus}\Osh_{D}(a) \stackrel{n}{\oplus} \Osh_{D}(a-1)$$ 
with $a\le 0$ and $m+n=(r-k)$.  

\medskip
\noindent
The curve $D$ is a $(-1)$ curve on the surface; and so setting 
$\Sbndle^{\prime}=f^{*}\Sbndle\otimes\Osh_{\Cfam}(aD)$ we have
$$ \Sbndle^{\prime}|_D = \,
\stackrel{m}{\oplus}\Osh_{D} \stackrel{n}{\oplus} \Osh_{D}(-1).$$ 
If $S_0$ is the bundle $f^{*}\Sbndle$ restricted to $C_0$, then 
$\Sbndle^{\prime}$ restricted to $C_0$ is generated by the sections of $S_0$
vanishing to order $(-a)$ at $p$. 

\medskip
\noindent
Shrink $B$ so that $b_0$ is the only singular fibre.
The bundle $\Sbndle^{\prime}$ then satisfies the conditions
of Lemma \ref{base_change_lemma} and therefore $g_{*}\Sbndle^{\prime}$ gives
a family of vector bundles over $C\times B$.  This bundle also has an inclusion
into $E$ induced by the map $\Sbndle^{\prime} = f^{*}\Sbndle\otimes \Osh(aD)
\rightarrow f^{*}\Sbndle$.

\medskip
\noindent
The quotient of $E$ by $g_{*}\Sbndle^{\prime}$ gives (using the same argument
as in Theorem \ref{(r-1)_map}) a flat family of quotients over $B$.  
Since the Hilbert
scheme is proper, the quotient at $b_0$ must be the correct one parameter limit,
and all that remains is to compute what this limit is.

\medskip
\noindent
To compare the restriction of $g_{*}\Sbndle^{\prime}$ to $C_{b_0}$ 
(on $C\times B$) and the restriction of $\Sbndle^{\prime}$ to $C_0$ 
(on the family $\Cfam$) we use the fact that formation of 
$g_{*}\Sbndle^{\prime}$ commutes with base change from $B$, and so we
only have to compute this for the contraction map
$$g_{b_0}: \Cfam_{b_0}=(C_0 \sqcup_{p} D) \longrightarrow C_0 = C_{b_0}.$$ 
Choose a local basis $s_1^{\prime},\ldots,s^{\prime}_{(r-k)}$ at $p$ of 
$\Sbndle^{\prime}$ restricted to $(C_0 \sqcup_p D)$, compatible with the 
splitting of $\Sbndle^{\prime}$ along $D$. Since the bundle $\Osh_{D}(-1)$ has 
no global sections, we see that $g_{*}\Sbndle^{\prime}$ is generated by
$\{s^{\prime}_1,\ldots,s^{\prime}_m,x\cdot s^{\prime}_{m+1},\ldots,
x\cdot s^{\prime}_{(r-k)}\}$.

\medskip
\noindent
The image of $\Sbndle^{\prime}|_{C_{b_0}}$ in $S_0$ is 
(as mentioned above) all the sections of $S_0$ vanishing to order $(-a)$ 
at $p$. The image of 
$g_{*}\Sbndle^{\prime}|_{C_0}$ is therefore $m$ sections vanishing to order
$(-a)$ and $n=(r-k-m)$ sections vanishing to order $(-a+1)$, as claimed in
the proposition \doneproof.

\begin{proposition}\label{one_(-2)_case}
Suppose that we have a family $\Cfam\longrightarrow B$ satisfying the same
conditions as in Proposition \ref{generic_splitting_case} but that this time
the restriction of $f^{*}\Sbndle$ to the rational curve splits as 
$(r-k-1)$ copies of $\Osh$ and one copy of $\Osh(-2)$, 
then the limiting quotient $\Qfam_{b_0}$ is also given by the rule above.
\end{proposition}
\noindent
This means that it is possible to choose a basis $\{s_1,\ldots,s_{(r-k)}\}$ for
$f^{*}\Sbndle$ restricted to $C_0$ (compatible with the splitting on the
rational curve) so that the limit is the quotient by 
$S^{\prime}=\{s_1,\ldots,s_{r-k-1}, x^2s_{(r-k)}\}.$ 

\medskip
\noindent
{\bf Remark.} In this case the subbundle $S^{\prime}$ {\em does} depend
on the basis chosen near $p$. For example, if $(r-k)=2$ the bases
$\{s_1,s_2\}$ and $\{s_1 + xs_2, s_2\}$ can both be compatible with the
splitting (at $p$) along the rational curve, but the subsheaves generated by
$\{s_1,x^2s_2\}$ and $\{s_1 + xs_2,x^2s_2\}$ are not the same.

\medskip
\noindent
{\em Proof.} Let $D$ be the rational curve and $C_0$ the section as before.  
The restriction of
$f^{*}\Sbndle$ to $D$ is
$$ f^{*}\Sbndle|_D =\, \stackrel{r-k-1}{\oplus}\Osh_D \oplus \Osh_D(-2).$$
Let $\Sbndle^{\prime}$ be the kernel (on $\Cfam$) of the projection onto the
$\Osh_D(-2)$ factor, so that we have
$$ 0\longrightarrow \Sbndle^{\prime} \longrightarrow f^{*}\Sbndle
\longrightarrow \Osh_D(-2) \longrightarrow 0.$$
We restrict this defining sequence to $D$ to get
$$ 0\longrightarrow \mathcal{T}or_1(\Osh_D(-2),\Osh_D) 
\longrightarrow \Sbndle^{\prime}|_D
\longrightarrow f^{*}\Sbndle|_D \longrightarrow \Osh_D(-2) \longrightarrow 0.$$
The map $f^{*}\Sbndle|_D \longrightarrow \Osh_D(-2)$ is 
just the projection we started with, and since 
$\mathcal{T}or_1(\Osh_D(-2),\Osh_D)=\Osh_D(-1)$, we have
$$ 0\longrightarrow \Osh_D(-1) \longrightarrow 
\Sbndle^{\prime}|_D
\longrightarrow \stackrel{r-k-1}{\oplus}\Osh_D \longrightarrow 0.$$
This sequence splits, and the splitting is canonical, since the trivial 
subbundle of $\Sbndle^{\prime}|_D$ is exactly the subbundle generated by global
sections.

\medskip
\noindent
Choose a basis $\{s^{\prime}_1,\ldots,s^{\prime}_{(r-k)}\}$ near $p$ for 
$\Sbndle^{\prime}$ restricted to $\Cfam_{b_0}$ compatible with the splitting
on $D$ (the basis element $s^{\prime}_{(r-k)}$ should be associated to the
$\Osh(-1)$ subbundle).

\medskip
\noindent
The image of $\Sbndle^{\prime}$ in $f^{*}\Sbndle$ restricted to $D$ is
$\stackrel{(r-k-1)}{\oplus} \Osh_D$.  This shows that we can pick a basis 
$\{s_1,\ldots,s_{(r-k)}\}$ for $f^{*}\Sbndle$ restricted to $C_0$ so that
the induced map $\Sbndle^{\prime}|_{C_0} \rightarrow S_0=f^{*}\Sbndle|_{C_0}$
is given by $s^{\prime}_i \mapsto s_i$ for $i=1,\ldots,(r-k-1)$ and
$s^{\prime}_{(r-k)} \mapsto x\cdot s_{(r-k)}$ for the remaining basis element.

\medskip
\noindent
We again pushforward the bundle $\Sbndle^{\prime}$ to get a bundle on 
$C\times B$ which we use to construct a flat family of quotients over $B$.
Following the same type of argument as in Proposition 
\ref{generic_splitting_case}, we see that $g_{*}\Sbndle^{\prime}|C_{b_0}$
is the subbundle of $\Sbndle^{\prime}|_{C_0}$ generated by
$\{s^{\prime}_1,\ldots,s^{\prime}_{(r-k-1)},x\cdot s^{\prime}_{(r-k)}\}$.

\medskip
\noindent
Using our previous description for the image of $\Sbndle^{\prime}|_{C_0}$ in 
$S_0=f^{*}\Sbndle|_{C_0}$ we see that the image $S^{\prime}$ of 
$g_{*}\Sbndle^{\prime}$ in $S_0$ is $\{s_1,\ldots,s_{(r-k-1)},x^2 s_{(r-k)}\}$
 \doneproof.

\medskip
\noindent
{\bf Remark.}  A combination of the arguments in Propositions 
\ref{generic_splitting_case} and \ref{one_(-2)_case} shows that the 
``rule'' for computing the limit also works if the subbundle splits as
$m_0$ copies of $\Osh_D(a)$,  $m_1$ copies of $\Osh_D(a-1)$, and 
$m_2$ copies of $\Osh_D(a-2)$.

\medskip
\noindent
The actual subbundle $S^{\prime}$ of $S_0$ produced 
depends on the bundle $\Sbndle^{\prime}$ which is in turn dependent on
the first order information of the curve $\Cfam_{b_0}$ in the family.
Different families {\em can} produce different limits, and this
is the reason that we cannot in general define a morphism from 
$\Mg(\Grk,\beta_d)$ to $\Quot$.  

\medskip
\noindent
The easiest way to see that the limit varies according to the family 
is to take a second limit consisting entirely of singular curves.  
If we stratify
the singular curves by the way the restriction of the universal subbundle 
decomposes on the rational tails,  then
the reason why we cannot extend the map is very clear:  the limit of the 
rules for computing the quotient on one strata is not the rule for computing
the quotient on a limiting strata.

\medskip
\noindent
The simplest case of this is to take a degeneration of the decomposition
$\Osh(-1) \oplus \Osh(-1)$ to the decomposition $\Osh \oplus \Osh(-2)$.

\begin{theorem}\label{dont_work_in_general}
If $k\ne (r-1)$ then in general there is no map from the space
$\Mg(\Grk,\beta_d)$
to ${\rm Quot}_{k,d}(E)$ extending the map on smooth curves.
\end{theorem}

\noindent
{\em Proof.} We show the case $(r-k)=2$; the idea works for all $k\ne(r-1)$.

\medskip
\noindent
Pick $d$ large so that by Theorem \ref{there_can_be_only_one} there is
a component of dimension $dr-ke-k(r-k)(g-1)$ 
whose generic point is a smooth point of the component, and corresponds
to a vector bundle quotient.

\medskip
\noindent
Let $C_0$ be a section corresponding to one of these points, 
and $S_0$ the restriction of the universal subbundle to $C_0$.
Pick a point
$p$ of $C_0$ and glue on a rational curve $D$.  Extend the inclusion
$C_0\hookrightarrow \Grk$ by mapping the curve in so that the universal
subbundle restricts to $D$ as $\Osh \oplus \Osh(-2)$. 

\medskip
\noindent
Because the point on $\Mg(\Grk,\beta_d)$ corresponding to $C_0$ was a 
smooth point of the component, a tangent space calculation shows that
this new stable map is a smooth point of
$\Mg(\Grk,\beta_{d+2})$ and lies on a component of dimension
$(d+2)r-ke-k(r-k)(g-1)$.  The dimension of the locus of nearby nodal curves is
one less than this, so the general point of this component corresponds to 
a smooth map. 

\medskip
\noindent
Take a smooth one parameter family $B_1$ containing our stable map as a 
fibre $b_0$.  Proposition \ref{one_(-2)_case} now shows that there is a basis
$\{s_1,s_2\}$ of $S_0$ near $p$ so that the limiting quotient is given by
dividing out by the subsheaf $S'_1=\{s_1,x^2s_2\}$.

\medskip
\noindent
Starting with the same stable map, we now construct another family $B_2$, by
fixing the section $C_0$, the point of attachment $p$ and varying only the
map of the rational curve $D$ into the Grassmannian fibre.  For a general 
variation, the decomposition of the universal subbundle will be $\Osh(-1)
\oplus \Osh(-1)$.  Let $b_0$ again stand for the fibre with the original
stable map.

\medskip
\noindent
Since all of the curves in our family can also be deformed to smooth curves,
Proposition \ref{generic_splitting_case} shows that for 
points where the decomposition is $\Osh(-1)\oplus\Osh(-1)$ the corresponding
point in $\Quot$ must be the quotient by $S^{\prime}_2=\{xs_1,xs_2\}$.  As
the point $b$ varies in $B_2$ the subbundle $S^{\prime}_2$ does not change,
and the inclusion $S^{\prime}_2\hookrightarrow
E$ is clearly the same away from $p$, and therefore the same at $p$ as well.
Therefore this is a {\em constant} family of quotients, and the limit
at $b_0$ in this family is {\em also} the quotient by $S^{\prime}_2$.

\medskip
\noindent
The two quotients are different, and this shows that (unless we blow up) 
we cannot extend the map on smooth curves to all of $\Mg(\Grk,\beta_{d+2})$ 
\doneproof.

\medskip
\noindent
{\bf Remark.} It is not clear if the fact that the splitting can change
is the only obstruction to extending the map.  In the two cases where the
decomposition {\em cannot} change (when $k=(r-1)$ or when $d=d_k$ or $d=d_k+1$) the map
can always be defined.

\medskip
\noindent
{\bf Remark.} For a general rational curve in a Grassmannian, the splitting
of the restriction of the universal subbundle is as ``democratic'' as 
possible -- no $a_i$ will differ by more than one from any other.  
Proposition \ref{generic_splitting_case} is therefore a reflection of the
fact that rational maps extend into codimension one subsets.

\bigskip
\noindent
\example \label{strange_limit}
Let $L$ be any line bundle, and $p$ any point on a curve $C$.  Set
$$ S=L\oplus L(p)\mbox{,\ \ \  and\ \ \ } E=L\oplus L(p)\oplus L(3p).$$
Let $b$ be a coordinate on the affine line, and define a family of
inclusions $S\stackrel{i_b}{\longrightarrow} E$ by
$$
\begin{array}{lcrcccl}
L & \stackrel{i_b}{\longrightarrow} & (&b^2\cdot \mbox{Id}_L,& \cdot \sigma_p,& 0&) \\
L(p) & \stackrel{i_b}{\longrightarrow} & (&0, & b\cdot\mbox{Id}_{L(p)}, &\cdot\sigma_p^2&)
\end{array}
$$

\medskip
\noindent
where $\sigma_p$ is the section of $\Osh_C(p)$ and ``$\cdot$'' indicates 
multiplication by the section.  If $x$ is a coordinate 
at $p$, then locally the inclusion is given by the matrix
$$\left(
\begin{array}{ccc}
b^2 &  0 \\
x & b \\
0 &  x^2
\end{array}
\right)$$
If $b\ne0$ then the inclusion is rank $2$ everywhere, and gives a vector bundle
quotient.  When $b=0$ the quotient has torsion, and looks like $E$ divided by
a subsheaf of the type $\{xs_1,x^2s_2\}$.

\medskip
\noindent
For $b\ne0$ we can consider this as a family of stable maps, and take the 
one parameter limit in $\Mg(\Grk,\beta_d)$. To do this, blow up the family
$C \times {\mathbb A}^1_b$ at $(p,0)$ and pull back both bundles to this
surface.  The inclusion of $S$ into $E$ now drops rank along
the exceptional divisor $D$, and if we take the saturation of $S$ in $E$
on the blow up, this gives the subbundle of the limiting stable map.
Computation in coordinates on the blowup shows that this subbundle 
restricted to $D$ decomposes as $\Osh_D \oplus \Osh_D(-3)$.

\medskip
\noindent
Since we know that the corresponding torsion quotient is {\em not} of 
the type where we divide out by something of the form $\{s_1,x^3s_2\}$, 
this shows that the decomposition of the subbundle on the rational tails 
has no relation in general to the corresponding limiting quotient in 
one parameter families.

\medskip
\noindent
{\bf Remark.} In the proof of Proposition \ref{one_(-2)_case} there was one
step where an exact sequence of bundles on $D$ split.  It is the failure of
this kind of sequence to split in general which prevents the rule from being true.
The example above is the simplest case where this splitting (and therefore
the rule) fails.

\medskip
\noindent
To finish, here is a result that properly belongs in the previous section,
but appears here because we need to use our picture of the relationship
between $\Mg(\Grk,\beta_d)$ and $\Quot$, especially pushing forward
vector bundles, and computing one parameter limits.

\medskip
\noindent
\begin{proposition}
\label{connected_Mg}
There is a $d_{\cal{M}}=d_{\cal{M}}(E,k)$ so that $\Mg(\Grk,\beta_d)$ is 
connected for all $d\ge d_{\cal{M}}$.
\end{proposition}

\medskip
\noindent
{\em Proof.} Choose $d_{\cal{M}}$ so that it is larger than the $d_{Q}$ of 
Theorem \ref{irred_of_Quot} and also large enough so that every component
of $\Mg(\Grk,\beta_d)$ has dimension at least $k(r-k)+1$.

\medskip
\noindent
For any $d\ge d_{\cal{M}}$, let $X$ be a component of $\Mg(\Grk,\beta_d)$.
By Lemma \ref{breaking_lemma} $X$ has a nonempty divisor $Y$ of reducible
maps.  Suppose that a point of a component of this divisor is
given by a section $C_0$ corresponding to a quotient of degree $(d-\delta)$
and an attached rational curve making up the difference in degree.  By
deforming the rational curve alone in the Grassmannian fibre, we can break
it into a tree of $\delta$ rational curves, each accounting for one in 
degree (this is a degeneration, and the resulting stable map is still a
point of $X$).

\medskip
\noindent
Now move the $\delta$ rational curves ``off'' the tree so that they are 
connected to $C_0$ at $\delta$ different points (it is at this stage that
we may be moving from one component to another).  The result is a section
$C_0$ with $\delta$ rational curves attached at $\left\{p_1,\ldots,
p_{\delta}\right\}$ on $C_0$, and so that each curve accounts
for one in degree.  Let $(C',f)$ be the resulting stable map into the
relative Grassmannian.  The claim is that this point is in the same connected
component as the locus parametrizing vector bundle quotients.

\medskip
\noindent
To see this, we use the fact that the pushforward $g_{*}(f^{*}\cal{\Fbndle})$ of
the pullback of the tautological quotient bundle on $\Grk$ gives a quotient
of $E$. This happens since (in this situation) the pullback of the canonical 
subbundle
$f^{*}\Sbndle$ decomposes as the sum of $\Osh$'s and one $\Osh(-1)$ on
each of the $\delta$ rational tails (see the discussion after Lemma 
\ref{base_change_lemma}).  The resulting pushforward $F'$ is a quotient with
torsion.  The torsion free part is the vector bundle quotient corresponding
to the section $C_0$, and the torsion part is supported at the points 
$p_1,\ldots, p_{\delta}$, and is of length one at each of those points.

\medskip
\noindent
Since $d_{\cal{M}}$ is larger than $d_Q$, $\Quot$ is irreducible, 
and so we can find a smooth one parameter family $B$ of quotients so that the 
general point is a vector bundle quotient and the special point $q\in B$ 
is the torsion quotient $F'$.  Since the general point of $B$ is a vector 
bundle quotient, we can also consider this a curve
in $\Mg(\Grk,\beta_d)$, and we take the one parameter limit approaching
$q\in B$.  To compute the limit we just need to blow up the family $C\times B$
at the points $\{(p_1,q),\ldots,(p_{\delta},q)\}$ since this is where the
quotient map drops rank, and take the saturation of the universal subbundle
pulled back to the blown up surface.  The limit is not necessarily the stable
map $(C',f)$ but it {\em is} a stable map with section $C_0$ (the {\em same}
$C_0$) and $\delta$ rational curves attached at $p_1,\ldots,p_{\delta}$ each
one accounting for one in degree.

\medskip
\noindent
This curve can be deformed (by moving the rational curves in the 
Grasmannian fibres) to $(C',f)$, and so we see that the component $X$ is in
the same connected component as the locus parametrizing vector bundle 
quotients.  Since $X$ was arbitrary, this means that $\Mg(\Grk,\beta_d)$ is
connected \doneproof.

\bigskip

\bigskip
\noindent
\mysection{Application to Base Point Freeness on Moduli Spaces}

\medskip
\noindent
In this section we assume that $C$ is a smooth projective complex 
curve of genus $g \geq 2$. Let 
$SU_{C}(r)$ denote the moduli space of semistable rank $r$ vector bundles on $C$ of trivial
determinant. The Picard group of $SU_{C}(r)$ is isomorphic to $\ZZ$ (see \cite{Drezet} Theorem $B$), 
with ample generator $\mathcal{L}$ (called the \emph{determinant}
bundle).  

\medskip
\noindent
For any
bundle $E$ in $SU_C(r)$, if there is a bundle $E'$ of rank $p$ and slope $(g-1)$
such that $h^0(E\otimes E')=0$, then $E$ is {\em not} a basepoint for 
$|\mathcal{L}^p|$. (It is not known if the converse is true -- this is the
content of the geometric formulation of the 
Strange Duality conjecture, see e.g. \cite{Beauville} \S8).

\medskip
\noindent
If we treat an element of $H^0(E\otimes E')$ as a homomorphism from $E^{*}$ to
$E'$, then we can factor the homomorphism as $ E^{*} \longrightarrow F 
\hookrightarrow E',$ where the map to $F$ is a surjection. 

\medskip
\noindent
Suppose that we fix a quotient bundle $F$ of rank $k$ and degree $d$, then 
we can ask for the dimension of stable bundles $E'$ of rank $p$ and 
slope $(g-1)$
which fit in the diagram
$$ 0\longrightarrow F \longrightarrow E' \longrightarrow Q \longrightarrow 0.$$
Counting dimensions of the quotient sheaf $Q$ and the extension 
space $\Ext^1(Q,F)$ shows that these bundles are of codimension at least $dp+1$ in 
the moduli space of bundles of rank $p$, slope $(g-1)$.
(In fact, for the purpose of dimension counts, $Q$ can be thought of a being 
a stable bundle -- for details and further references related to all these facts,
please see \cite{Popa}, especially \S4).

\medskip
\noindent
This gives the following lemma on basepoint freeness:

\begin{lemma}
\label{bpf_lemma}
For a bundle $E$ in $SU_C(r)$, if $p$ is a number so that 
$$\dim \left\{{\mbox{
\begin{minipage}{5cm}
\noindent
Moduli of rank $k$, degree $d$, quotient bundles $F$ of $E$
\end{minipage}}}
\,\right\} \le dp$$
for all $k = 1,\ldots,(r-1)$, and all $d$, then $E$ is not a basepoint of
$|\mathcal{L}^p|$.
\end{lemma}

\noindent
The dimension above is the dimension of the moduli of the quotient bundles
$F$, throwing away the data of the surjections $E\rightarrow F$.  The space
of quotients (i.e., this time including the data of the surjections) provides 
an upper bound for the dimension of this moduli.  We
have $d_k(E)\ge 1$ for all $k$ since $E$ is stable of degree $0$, 
and using our upper bound \ref{Quot_dim} gives the following (note that the
argument works literally only for $r\geq 3$, but it is well known that 
$\mathcal{L}$ is already very ample on $SU_{C}(2)$):

\begin{theorem}\label{3}
The series $|\mathcal{L}^{p}|$ on $SU_{C}(r)$ 
is base point free for $p\geq [\frac{r^{2}}{4}]$.
\end{theorem}

\medskip
\noindent
{\bf Remark.}
This is a small improvement on a similar bound in \cite{Popa}.  In the same
paper, the following question appeared:

\medskip
\noindent
{\bf Question.}\label{mihneas_conjecture}
Is the series $|\mathcal{L}^p|$ base point free on 
$SU_C(r)$ for $p\ge (r-1)$?

\medskip
\noindent
Even though we have examples of bundles with large families of quotients,
we know of no potential counterexample (i.e. modulo Strange Duality) 
to this conjecture. There are two difficulties 
in constructing good examples in this situation. 
The first is that we would like the family of quotients to have maximal 
(or at least large) variation in moduli. The second (and more difficult) is that
the dimension count in Lemma \ref{bpf_lemma} does not take into account the
fact that different extensions involving different $Q$'s may produce the
same stable bundle $E'$.  This contribution is difficult to estimate.

\medskip
\noindent
For instance, the bundles $G_0$ constructed in Example 
5.\ref{extension_example} are stable bundles of rank $r$ and trivial 
determinant with at least a $(k-1)(r-k-1)$ dimensional family of rank $k$,
degree $1$ quotient bundles.  As stated in
the example, these families of quotients vary maximally in moduli, and
so seem as if they should provide a potential 
counterexample to the conjecture.
However, all the quotients $F$ are extensions of the form
$$0\rightarrow \stackrel{k-1}{\oplus} {\cal O} \rightarrow F \rightarrow L
\rightarrow 0,$$
where $L$ is a fixed line bundle of degree $1$.  
The general stable bundle $E'$ of slope $(g-1)$ has no global sections, 
and also will not accept a nonzero map from $L$, and therefore any morphism
of the type $F\rightarrow E'$ is the zero morphism.

\medskip
\noindent
In any construction we can make, either the dimension of the quotients is
too small, or the quotient bundles produced are too special.  The vague
principle seems to be that quotients of a fixed bundle cannot be too general,
and the Question above is a quantitative expression of this.

\medskip
\bigskip
\noindent
The improved upper bound also gives a few new cases of a conjecture 
made in \cite{Popa} 5.5 about linear 
series on the moduli space $U_{C}(r,0)$ of semistable vector bundles of 
rank $r$ and degree $0$ on $C$.

\medskip
\noindent
Once a line bundle $N\in {\rm Pic}^{g-1}(C)$ is fixed, we can define (see \cite{Drezet} 7.4.2)
a \emph{generalized theta divisor} $\Theta_{N}$ on $U_{C}(r,0)$, given set-theoretically by:
$$\Theta_{N}=\{E~|~h^{0}(E\otimes N)\neq 0\}.$$
The conjecture mentioned above asserts that the linear series $|k\Theta_{N}|$ should be base point free 
for $k\geq r+1$. In the case $r=1$, this is the classical statement that $|2\Theta_{N}|$ is base 
point free on the Jacobian of $C$, and in \cite{Popa} 5.4 a proof is also given for $r=2$ and $3$.
By arguments in \cite{Popa2} 5.3, $|k\Theta_{N}|$ is base point free as long 
as $k\geq r+1$ and $|\mathcal{L}^{k}|$ is base point free on $SU_{C}(r)$. Concluding,  
Theorem \ref{3} provides a solution to the conjecture up to the case of rank $5$ vector bundles:

\begin{corollary}
$|k\Theta_{N}|$ is base point free on $U_{C}(r,0)$ for $k\geq r+1$ if $r\leq 5$.
\end{corollary}

\medskip
\noindent
Department of Mathematics, Harvard University, \\
One Oxford Street, Cambridge, MA 02138, U.S.A
\newline
email: mpopa@math.harvard.edu

\medskip
\noindent
Department of Mathematics, University of Michigan, \\
525 East University, Ann Arbor, MI 48109-1109, U.S.A.
\newline
email: mikeroth@umich.edu

\end{document}